\documentclass[11pt,twoside]{article}

\usepackage{authblk}
\usepackage{fancyhdr}
\usepackage{amsfonts,epsfig,graphicx}
\usepackage{afterpage}
\usepackage{amsmath,amssymb,amsthm} 
\usepackage{fullpage}
\usepackage[T1]{fontenc} 
\usepackage{epsf} 
\usepackage{graphics} 
\usepackage{amsfonts,amsmath}
\usepackage[]{natbib} 
\usepackage{psfrag,xspace}
\usepackage{color,etoolbox}
\usepackage[skip=1pt]{caption}
\usepackage{subfigure} % subfigures
\usepackage{hyperref}[] % hyperlink
\hypersetup{
	colorlinks = true,
	linkcolor = blue, %Colour of internal links
	filecolor = magenta,
	urlcolor = blue, %Colour for external hyperlinks
	citecolor = blue %Colour of citations
}

\newcommand{\mP}{\mathbb{P}}
\newcommand{\mE}{\mathbb{E}}

\newcommand{\mS}{\mathcal{S}}
\newcommand{\bX}{\textbf{X}}

\newcommand{\convP}{\overset{p}{\longrightarrow}}
\newcommand{\convD}{\overset{d}{\longrightarrow}}

\newtheorem{theorem}{Theorem}[section]
\newtheorem{lemma}{Lemma}[section]

\newtheorem{corollary}{Corollary}[section]
\newtheorem{remark}{Remark}[section]
\newtheorem{example}{Example}[section]

\begin{document}

\title{\textbf{\Large Multinomial Goodness-of-Fit Based on $U$-Statistics: High-Dimensional Asymptotic and Minimax Optimality}}

\author[1]{Ilmun Kim}
\affil[1]{Department of Statistics \& Data Science, \protect\\ Carnegie Mellon University, Pittsburgh, PA 15213, USA \protect\\ {\normalsize\href{mailto:ilmunk@andrew.cmu.edu}{ilmunk@andrew.cmu.edu}}}
\date{}

\maketitle

\begin{abstract}
	We consider multinomial goodness-of-fit tests in the high-dimensional regime where the number of bins increases with the sample size. In this regime, Pearson's chi-squared test can suffer from low power due to the substantial bias as well as high variance of its statistic. To resolve these issues, we introduce a family of $U$-statistic for multinomial goodness-of-fit and study their asymptotic behaviors in high-dimensions. Specifically, we establish conditions under which the considered $U$-statistic is asymptotically Poisson or Gaussian, and investigate its power function under each asymptotic regime. Furthermore, we introduce a class of weights for the $U$-statistic that results in minimax rate optimal tests.
\end{abstract}

\vskip 1em

\begin{small}
\end{small}

\vskip 1em

\section{Introduction}
Suppose that there are $n$ independent random vectors $\bX_1 = (X_{1,1}, \ldots , X_{1,d}), \ldots , \bX_n=(X_{n,1}, \ldots , X_{n,d})$  from a multinomial distribution with unknown parameters $\pi = (\pi_{1}, \ldots , \pi_d) \\ \in \Omega$ and
\begin{align*}
\Omega = \Big\{ (\pi_1, \ldots , \pi_d) \in [0,1]^d : \sum_{j=1}^d \pi_j = 1\Big\}.
\end{align*}
Given a specific choice of parameter vector $\pi_0 = (\pi_{0,1}, \ldots, \pi_{0,d}) \in \Omega$, the goodness-of-fit test for multinomial distributions is concerned with distinguishing
\begin{align} \label{Eq: Hypothesis}
H_0: \pi_0 =\pi \quad  \text{versus} \quad H_1: \pi_{0} \neq \pi.
\end{align}
Pearson's chi-squared statistic \citep{pearson1900x} is one of the well-known test statistics for this problem. Let $Y_j = \sum_{i=1}^n I(X_{i,j}=1)$ for $j=1,\ldots, d$. Then Pearson's chi-squared statistic is defined by
\begin{align*}
\chi^2_n = \sum_{j=1}^d \frac{\left( Y_j - n\pi_{0,j} \right)^2}{n \pi_{0,j}}.
\end{align*}
The properties of $\chi_n^2$ have been well-studied in a classical low-dimensional setting \citep{lehmann2006testing, read2012goodness,balakrishnan2013chi}. For instance, the test based on $\chi_n^2$ is asymptotically optimal against local alternatives when $d$ is fixed \citep[see, e.g. Chapter 14 of][]{lehmann2006testing}. However, in the high-dimensional regime where the dimension is comparable with or much larger than the sample size, $\chi_n^2$ suffers from the fact that it can have substantial bias for the testing problem. In other words, the power of the test can be much smaller than the significance level $\alpha$ against certain local alternatives. The major cause of the testing bias is due to the expected value of $\chi_n^2$:
\begin{align*}
\mE\left[ \chi_n^2 \right] = d - 1 + \sum_{j=1}^d  \frac{\pi_j - \pi_{0,j}}{\pi_{0,j}} + \frac{n-1}{n} \sum_{j=1}^d \frac{(\pi_j - \pi_{0,j})^2}{\pi_{0,j}}.
\end{align*}
When the null is not uniform, it is possible to observe $\mE_{H_1}[\chi_n^2]<\mE_{H_0} [\chi_n^2]$, which can trigger a significant bias problem of $\chi_n^2$ for some $\alpha$ level. This bias problem becomes more serious when the dimension is large but the sample size is small \citep[see,][for details]{haberman1988warning}.

\begin{figure}[!t]
	\center
	\includegraphics[width=\textwidth]{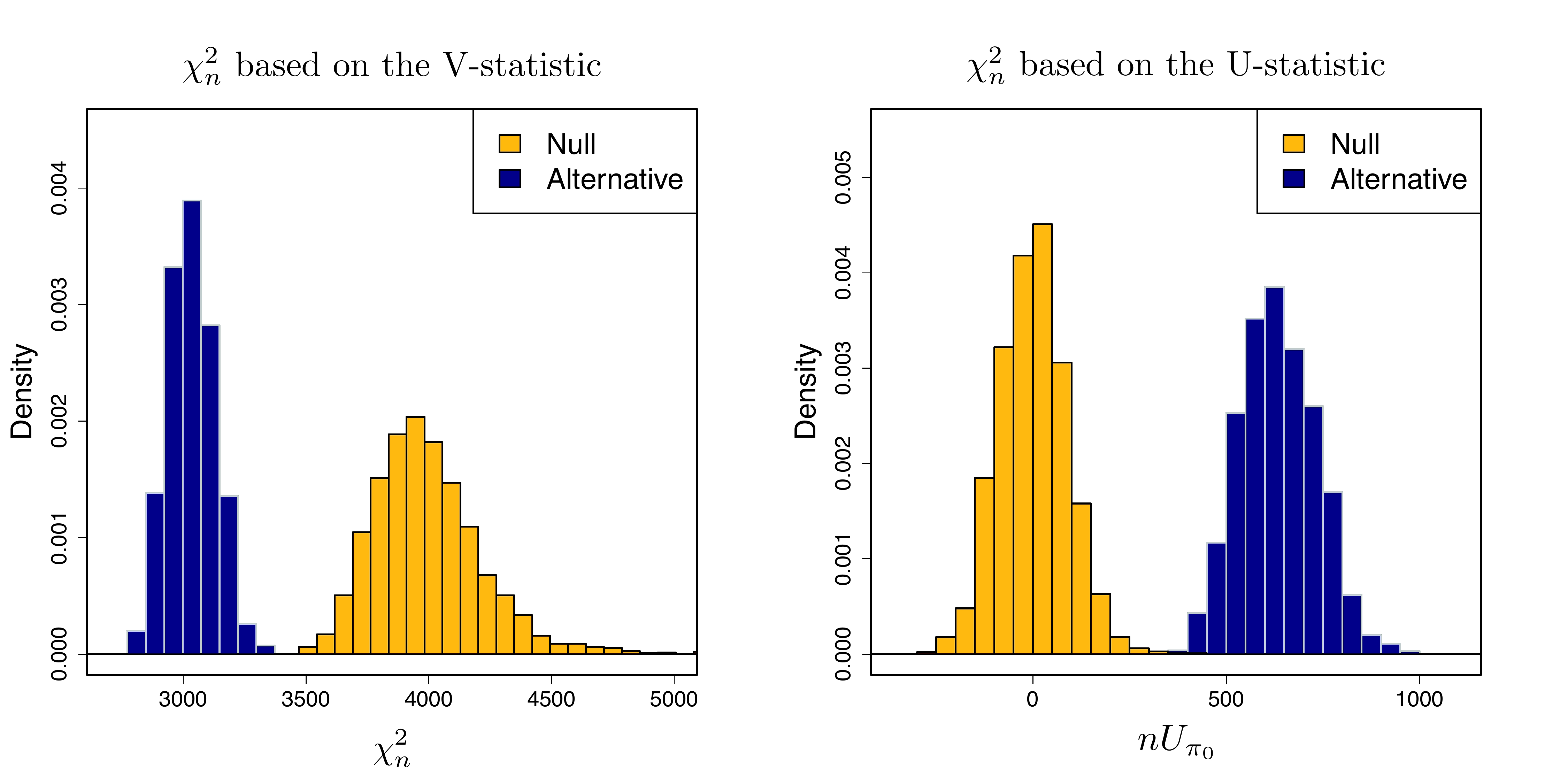}
	\caption{An illustration of the bias issue of Pearson's $\chi_n^2$ test in the high-dimensional regime. For the simulation, the null and the alternative are chosen $\pi_{0,i} \propto i$ and $\pi_i \propto i^5$, respectively. We take the sample size $n=800$ and the dimension $d=4000$. Since the null is rejected when the test statistic is greater than a certain quantile of the null distribution, we see from the left panel that the $\chi_n^2$ test is substantially biased. On the other hand, the right panel shows that the modified $\chi_n^2$ based on the $U$-statistic can have significant power in this example.} \label{Fig: V- and $U$-statistic}
\end{figure}

To avoid the testing bias caused by the expected value, we view Pearson's chi-squared statistic as a $V$-statistic (Lemma~\ref{Lemma: Pearson V-stat}) and consider a modified $\chi^2_n$ based on the $U$-statistic. From the basic property of $U$-statistics, the modified $\chi^2_n$ is an unbiased estimator of $\sum_{j=1}^d {(\pi_j - \pi_{0,j})^2}/{\pi_{0,j}}$ and its expectation becomes zero if and only if the null is true. As a result, the modified $\chi_n^2$ can have significant power in the high-dimensional regime where classical $\chi^2_n$ is substantially biased (Figure~\ref{Fig: V- and $U$-statistic}). 

Another limitation of $\chi_n^2$ in sparse multinomial settings is that it puts too much weight on small entries in $\pi_0$, and these small entries make the statistic highly unstable \citep{marriott2015geometry,valiant2017automatic,balakrishnan2017hypothesis}. In this case, one might need to consider different weights to obtain higher power of the test. Motivated by these observations, we consider a family of $U$-statistics of $||A^{1/2}(\pi - \pi_0)||_2^2$ where $A$ is some positive definite matrix.

The primary objective of this work is to investigate the limiting behavior of the proposed $U$-statistic in high-dimensions and determine a sufficient condition for $A$ under which the resulting test is minimax rate optimal for multinomial goodness-of-fit. 

\vskip .8em

\noindent \textbf{Main results.} The main results of this paper are as follows:
\begin{enumerate}
	\item \emph{Poissonian Asymptotic for the $U$-statistic (Section \ref{Sec: Poisson Limit})}: We establish conditions under which the $U$-statistic has a Poisson limiting distribution.
	\item \emph{Gaussian Asymptotic for the $U$-statistic (Section \ref{Sec: Global Gaussian Limits})}: We also establish conditions under which the $U$-statistic has a Gaussian limiting distribution.
	\item \emph{Global Minimax Optimality of the $U$-statistic (Section \ref{Sec: Minimax Upper Bound})}: We present a class of weight matrices $A$ resulting in the minimax optimal test based on the $U$-statistic.  
\end{enumerate}

\vskip .8em

\noindent \textbf{Related work.} A considerable amount of literature has been published on the high-dimensional behavior of $\chi_n^2$ \citep[e.g.][and the references therein]{tumanyan1954asymptotic,tumanyan1956asymptotic,steck1957limit, holst1972asymptotic,morris1975central,read2012goodness,rempala2016double}. Our work is especially motivated by \cite{rempala2016double} who present conditions of the Poissonian and Gaussian asymptotics for  $\chi_n^2$. One can generalize their result to our $U$-statistic framework by using the relationship between $U$- and $V$-statistics. However, their analysis is restricted to the case of the null hypothesis and does not easily generalize to other cases with different weights. Hence, we take different approaches to overcome such shortcomings. The present study is also closely related to the work by \cite{zelterman1986log,zelterman1987goodness} who proposes a modified $\chi_n^2$ to handle the testing bias of the chi-squared test. The modified statistic is given by 
\begin{align} \label{Eq: zelterman}
\phi_n = \chi^2_n  - \sum_{j=1}^d \frac{Y_j}{n \pi_{0,j}}.
\end{align} 
It can be shown that $\phi_n$ is equivalent to the proposed $U$-statistic up to some constant factors when we take the weight matrix as $A = \text{diag}(\pi_{0,1}^{-1}, \ldots, \pi_{0,d}^{-1})$ (Remark~\ref{Remark: Zelterman}), and thus $\phi_n$ falls into our $U$-statistic framework. \cite{diakonikolas2016collision} show that the collision-based test is minimax rate optimal for multinomial uniformity testing where $\pi_0 = (1/d,\ldots, 1/d)^\top$ (Remark~\ref{Remark: Collision-based test}). Their test statistic is a special case of our $U$-statistics by taking the identity weight matrix given in (\ref{Eq: U-statistic with the identity matrix}). We generalize their minimax result to an arbitrary null probability by providing a class of $A$ that leads to the minimax optimal test. Our result includes the truncated weight considered in \cite{balakrishnan2017hypothesis} as an example.

\vskip .8em

\noindent \textbf{Outline.}  The rest of the paper is organized as follows. In Section~\ref{Sec: Pearson based on $U$-statistic}, we view Pearson's chi-squared statistic as a $V$-statistic and provide a modified and generalized $\chi_n^2$ based on the $U$-statistic. In Section~\ref{Sec: Asymptotic Distribution in High-Dim}, we investigate the Poissonian and Gaussian asymptotics for the proposed $U$-statistic in the high-dimensional regime. In Section~\ref{Sec: Minimax Upper Bound}, we present a sufficient condition for $A$ that results in the minimax optimal test based on the $U$-statistic. In Section~\ref{Sec: Simulations}, we provide simulation studies. We summarize the results and discuss future work in Section~\ref{Sec: Conclusions}. Finally, all of the proofs and additional results are presented in Appendix~\ref{Sec: Appendix}.

\vskip 2em

\section{Pearson's Chi-squared Statistic based on the $U$-statistic} \label{Sec: Pearson based on $U$-statistic}
As mentioned earlier, when the null is not uniform, Pearson's chi-squared statistic can have $\mE_{H_1}[\chi_n^2] < \mE_{H_0}[\chi_n^2]$. This phenomenon can result in serious testing bias especially in the high-dimensional regime. To remove the main testing bias due to its expected value, we first view $\chi_n^2$ as a $V$-statistic and suggest the modified $\chi_n^2$ based on a $U$-statistic. To begin, consider the following second order kernel:
\begin{align} \label{Eq: chi-square kernel}
h_{\pi_0}(\bX_i, \bX_j)  = (\bX_i - \pi_0)^\top D_{\pi_0}^{-1} (\bX_j - \pi_0),
\end{align}
where $D_{\pi_0}$ is a $d \times d$ diagonal matrix with the diagonal entries $\{ \pi_{0,1}, \ldots, \pi_{0,d} \}$. Given $h_{\pi_0}$, we define the $V$-statistic as 
\begin{align*}
V_{\pi_0} = \frac{1}{n^2} \sum_{i=1}^n \sum_{j=1}^n h_{\pi_0}(\bX_i, \bX_j).
\end{align*}
In this setting, the next lemma shows that Pearson's chi-squared statistic is equivalent to the $V$-statistic defined with kernel $h_{\pi_0}$.

\begin{lemma} \label{Lemma: Pearson V-stat}
Pearson's chi-squared statistic has another representation as	
\begin{align*}
\chi^2_n = \sum_{j=1}^d \frac{\left( Y_j - n\pi_{0,j} \right)^2}{n \pi_{0,j}} = \frac{1}{n} \sum_{i=1}^n \sum_{j=1}^n h_{\pi_0}(\emph{\bX}_i, \emph{\bX}_j) = n V_{\pi_0}.
\end{align*}
\end{lemma}

\vskip .8em

As is well-known, a $V$-statistic is typically biased for estimating the population quantity. In order to remove the estimation bias of $V_{\pi_0}$, we consider a $U$-statistic defined as
\begin{align} \label{Eq: Ustat of chi-sq}
U_{\pi_{0}} = \binom{n}{2}^{-1} \sum_{1 \leq i < j \leq n} h_{\pi_0}(\bX_i, \bX_j).
\end{align}
It can be easily seen that the expected value of $U_{\pi_0}$ is always non-negative and equal to zero if and only if $\pi = \pi_0$. As a result, $U_{\pi_{0}}$ does not suffer from the testing bias arising from the expected value.

\begin{remark} \label{Remark: Zelterman}
	$U_{\pi_{0}}$ is closely related to the test statistic proposed by \cite{zelterman1986log,zelterman1987goodness}. In view of Lemma~\ref{Lemma: Pearson V-stat}, it is straightforward to show that these two statistics have the identity $\phi_n  = (n-1)(U_{\pi_0} - 1 )$; thus they have the exact same power. Unfortunately, even if $U_{\pi_{0}}$ does not have the problem of the expectation, it can still have the testing bias against certain alternatives. For instance, \cite{haberman1988warning} provides a case where $\phi_n$'s power is less than the significance level, which implies that $U_{\pi_0}$ also has the testing bias in the same case.
\end{remark}

There are some theoretical and empirical evidence to suggest that the scaling factor of $\chi_n^2$ might not be optimal in the high-dimensional setting \citep{marriott2015geometry,valiant2017automatic,balakrishnan2017hypothesis}. For example, when $\pi_0$ is highly skewed, $\chi_n^2$ can perform poorly as it is dominated by small domain entries of $\pi_0$. Therefore, one might need to consider different weights for $\chi_n^2$ to obtain higher power of the test in high-dimensions. In this context, we consider a family of $U$-statistics by considering a general weight matrix. The test statistic of our interest is given as
\begin{align} \label{Eq: U_A definition}
U_A = \binom{n}{2}^{-1}\sum_{1 \leq i < j \leq n} h_A(\bX_i, \bX_j),
\end{align}
where $h_A(\bX_i, \bX_j) = (\bX_i - \pi_0)^\top A(\bX_j - \pi_0)$ and $A$ is some positive definite matrix. In the next section, we study the limiting behavior of $U_A$ under different high-dimensional regimes.

\vskip 2em

\section{High-Dimensional Asymptotics} \label{Sec: Asymptotic Distribution in High-Dim}
The asymptotic behavior of Pearson's chi-squared statistic has been well studied in the literature \citep[][for a review]{ read2012goodness}. In the high-dimensional case where the dimension (the number of bins) increases with the sample size, \cite{rempala2016double} investigate both the Poisson and Gaussian approximations for $\chi_n^2$. Specifically, when $\pi_0$ is uniformly distributed, they show that $\chi_n^2$ is asymptotically Poisson when $n/\sqrt{d} \rightarrow c \in (0,\infty)$ and asymptotically Gaussian when $n / \sqrt{d} \rightarrow \infty$. One can use their results to establish similar asymptotics for the $U$-statistic in view of Lemma~\ref{Lemma: Pearson V-stat}. However, their analysis is restricted to the case of the null hypothesis.

In this section, we derive both null and alternative distributions of the considered $U$-statistic and present conditions for its high-dimensional limiting behavior. Note that, under the low-dimensional regime where $d$ is fixed, it is rather straightforward to obtain the limiting distribution of the considered $U$-statistic. For example, $U_A$ is the $U$-statistic, which has degeneracy of order one at the null hypothesis; thereby it converges to a weighted sum of independent chi-squared random variables with one degree of freedom \citep[see, e.g.,][for asymptotic results of $U$-statistics]{lee1990u}. More interesting and challenging might be the high-dimensional case where $U_A$ can have a Poisson or Gaussian limiting distribution, which will be studied in the following subsections.

\subsection{Poisson Approximation} \label{Sec: Poisson Limit}
We start with establishing conditions under which the $U$-statistic has an asymptotic Poisson distribution. It is worth noting that, since a Poisson random variable is supported on the non-negative integers, an arbitrary choice of $A$ does not necessarily yield a Poisson approximation for $U_A$ (even after being properly centered and scaled). For this reason, we focus on the simple case where $A$ is the identity matrix, i.e.
\begin{align} \label{Eq: U-statistic with the identity matrix}
U_I = \binom{n}{2}^{-1} \sum_{1\leq i < j \leq n} (\bX_i - \pi_0)^\top (\bX_j - \pi_0),
\end{align}
and study its asymptotic behavior. We briefly discuss the generalization of the identity matrix to an arbitrary matrix $A$ resulting in the Poisson asymptotic in Remark~\ref{Remark: Poisson Generalization}.

Let us start by defining some conditions which hold as $n,d \rightarrow \infty$ simultaneously:

\begin{itemize}\setlength{\itemindent}{.2in}
\item [\textbf{(P.1)}] $n^3 \big\{  \sum_{i=1}^d \pi_i^3 + ( \pi^\top \pi )^2 \big\} \rightarrow 0$. 
\item [\textbf{(P.2)}] $\binom{n}{2} \pi^\top \pi \rightarrow \eta_1$, $\binom{n}{2} \pi_0^\top \pi_0 \rightarrow \eta_0$ and $\binom{n}{2}  \pi^\top \pi_0 \rightarrow \eta_2$ where $\eta_i \in (0, \infty)$ for $i=0,1,2$.
\item [\textbf{(P.3)}] $n^3 \big\{ \sum_{i=1}^d \pi_i \pi_{0,i}^2   - \big(\sum_{i=1}^d \pi_i \pi_{0,i} \big)^2 \big\} \rightarrow 0$.
\end{itemize}

\vskip .8em

Let us consider the following decomposition of $U_I$:
\begin{align*}
\binom{n}{2} U_I = W - (n-1)\sum_{i=1}^n \left( \bX_i^\top \pi_0 + \pi_0^\top \pi_0\right),
\end{align*}
where $W= \sum_{1 \leq i < j \leq n} \bX_i^\top \bX_j$. Based on this decomposition, first note that condition (P.1) together with the first condition in (P.2) is to ensure that $W$ is approximately Poisson with mean $\eta_1$. The last two conditions in (P.2) as well as (P.3) are to guarantee that the remainder of $\binom{n}{2}U_I$ other than $W$ converges to $\eta_0 - 2 \eta_2$ in probability. Specifically, (P.3) is a sufficient condition under which the variance of $(n-1)\sum_{i=1}^n \bX_i^\top \pi_0$ converges to zero so that the asymptotic behavior of $\binom{n}{2}U_I$ is dominated by $W$.

Under the above conditions, we depict the limiting behavior of $U_I$ as follows:

\begin{theorem}[Poisson limiting distributions] \label{Thm: Poisson approximation}
Under the conditions \emph{{(P.1)}}, \emph{{(P.2)}} and \emph{{(P.3)}}, $U_I$ has a Poisson limiting distribution as 
\begin{align*}
\binom{n}{2} U_I \convD \text{\emph{Pois}} (\eta_1) - 2 \eta_2 + \eta_0.
\end{align*}
\end{theorem}

\vskip 0.8em

In the following corollaries, we demonstrate the above result under the uniform null and the piecewise uniform alternatives.

\vskip 0.8em

\begin{corollary}[Uniform null distribution] \label{Corollary: uniform null}
Suppose that we are under the uniform null, i.e. $\pi = \pi_{0} = (1/d, \ldots, 1/d)^\top$. If $n / \sqrt{d} \rightarrow \sqrt{2\eta_0} \in (0, \infty)$, then
\begin{align*}
\binom{n}{2} U_I \convD \text{\emph{Pois}} (\eta_0)  -  \eta_0.
\end{align*}
If $\eta_0 = 0$, then it converges to zero in distribution. 
\end{corollary}

\vskip 0.8em

\begin{corollary}[Piecewise uniform alternatives] \label{Corollary: Piecewise uniform alternative}
	Suppose that the null distribution is uniformly distributed. Consider $\omega_1, \omega_2 >0$ such that $\omega_1 + \omega_2 = 1$. For simplicity, assume that $d$ is even number (otherwise, let $\pi_{1,d}=0$) and consider the alternative distribution defined by
	\begin{align*}
	\pi_{1}  =  ( \underbrace{\omega_1/d, \ldots, \omega_1/d}_{\text{$d/2$ elements}}, \underbrace{\omega_2/d, \ldots, \omega_2/d}_{\text{$d/2$ elements}}).
	\end{align*}
	If $n / \sqrt{d} \rightarrow \sqrt{2\eta_0} \in (0, \infty)$, then under the alternative hypothesis, 
	\begin{align*}
	\binom{n}{2} U_I \convD \text{\emph{Pois}} (\eta_1) -  \eta_0,
	\end{align*}
	where $\eta_1 = \eta_0 (\omega_1^2 + \omega_2^2)/2$. 
\end{corollary}

\vskip .8em

From the above results, let us describe the asymptotic power function of $U_I$ under the Poissonian asymptotic. We assume that the null distribution is uniform where $n/ \sqrt{d} \rightarrow c \in (0, \infty)$; thereby the distribution of $\binom{n}{2}(U_I + 1/d)$ is approximated by the Poisson distribution. Let $c_\alpha \in \mathbb{Z}^+$ be a critical value such that 
\begin{align*}
\mP_{H_0}\left(\binom{n}{2}(U_I + 1/d) > c_\alpha \right) \leq \alpha,
\end{align*}
under the null. Then the power function of $U_I$ can be approximated by
\begin{align} \label{Eq: Power function of Poisson}
\beta_{n,d}(\pi) & =  \mP_{H_1} \left( \binom{n}{2} (U_I + 1/d)  > c_\alpha \right)    = \int_0^{2\eta{_1}} \frac{1}{\Gamma(c_\alpha+1)} y^{c_\alpha} \exp\left({-\frac{y}{2}} \right) dy + o(1),
\end{align}
against the alternatives that satisfy (P.1), (P.2) and (P.3). %Note that the integrand of (\ref{Eq: Power function of Poisson}) is always positive; therefore, under the Poisson regime, $U_I$ is asymptotically powerful as long as $\eta_1 > \eta_0$ for the uniform null distribution.

\begin{remark} \label{Remark: Poisson Generalization}
	The Poisson approximation for $U_I$ can be extended to a general statistic $U_A$ when the weight matrix $A$ is asymptotically close to $\sigma I$ for some $\sigma > 0$. Suppose that we are under the null hypothesis and the conditions (P.1), (P.2) and (P.3) are satisfied. For $\Sigma_0 = \text{\emph{diag}}(\pi_0) - \pi_0 \pi_0^\top$ and $D_\sigma = A - \sigma I$, assume that $n^2 \text{\emph{tr}} \{ (D_\sigma \Sigma)^2 \}
	\rightarrow 0$ as $n,d \rightarrow \infty$. Then the following holds by Chebyshev's inequality:
	\begin{align*}
	\binom{n}{2} \left( U_A - \sigma U_I \right) \convP 0 \quad  \text{and} \quad \binom{n}{2} U_A \convD \sigma \text{\emph{Pois}} (\eta_0) - \sigma \eta_0.
	\end{align*}
\end{remark}

\subsection{Gaussian Approximation} \label{Sec: Global Gaussian Limits}
In this section, we study the asymptotic normality of $U_A$.  Without loss of generality, we further assume that $A$ is symmetric. Under the null hypothesis, the next theorem provides a sufficient condition under which $U_A$ is asymptotically Gaussian.

\begin{theorem}[Asymptotic normality of $U_A$ under the null] \label{Thm: Global Gaussian}
Suppose
\begin{align} \label{Eq: Gaussian Martingale Assumptions}
\frac{\emph{tr} ((A\Sigma)^4)}{ [\emph{tr} \{ (A\Sigma)^2 \} ]^2 } \rightarrow 0 \quad \text{and} \quad \frac{ \mE \big[ \big\{ h_A(\bX_1, \bX_2) \big\}^4 \big] + n\mE \big[ \big\{  h_A(\bX_1, \bX_2)  \big\}^2 \big\{ h_A(\bX_1, \bX_3)  \big\}^2  \big]}{ n^2 [\emph{tr} \{ (A\Sigma)^2 \} ]^2} \rightarrow 0.
\end{align}
Then, under the null hypothesis, we have
\begin{align*}
\sqrt{\binom{n}{2}}\frac{U_A}{\sqrt{\emph{tr} \big\{ (A\Sigma)^2 \big\} }}   \convD  \mathcal{N}(0,1).
\end{align*}
\end{theorem}

\vskip .8em

Recall that, in Corollary~\ref{Corollary: uniform null}, we established the limiting behavior of $U_I$ under the uniform null when $n/\sqrt{d} \rightarrow c \in [0,\infty)$. In the next corollary, we study the asymptotic normality of $U_I$ when $n/\sqrt{d} \rightarrow \infty$ under the uniform null.

\begin{corollary}[Uniform null distribution] \label{Corollary: Uniform Null Global Gaussian}
Suppose we are under the null hypothesis where $\pi_0$ is uniform. If $n/\sqrt{d} \rightarrow \infty$, then we have
\begin{align*}
\sqrt{\binom{n}{2}}\frac{U_I}{\sqrt{1/d(1-1/d)}}   \convD  \mathcal{N}(0,1).
\end{align*}
\end{corollary}

\vskip .8em

Let us now turn our attention to the alternative distribution of $U_A$ in the high-dimensional asymptotic. As in \cite{chen2010two}, we consider the following two scenarios under the alternative hypothesis:
\begin{itemize}\setlength{\itemindent}{.2in}
\item [\textbf{(S.1)}] (Strong Signal-to-Noise) $n^{-1}\text{tr}( (A\Sigma)^2)=  o\big( (\pi-\pi_0)^\top A \Sigma A(\pi - \pi_0) \big)$.
\item [\textbf{(S.2)}] (Weak Signal-to-Noise) $(\pi-\pi_0)^\top A \Sigma A(\pi - \pi_0) = o \big( n^{-1}\text{tr}( (A\Sigma)^2) \big)$.
\end{itemize}
To appreciate the given scenarios, let us decompose $U_A = U_{\text{quad}} + U_{\text{linear}}$ where
\begin{align*}
& U_{\text{quad}} =  \binom{n}{2}^{-1} \sum_{1\leq i < j \leq n} (\bX_i- \pi)^\top A (\bX_j - \pi),   \\[.5em]
& U_{\text{linear}} = \binom{n}{2}^{-1}  \sum_{1\leq i < j \leq n} \Big\{ (\bX_i - \pi)^\top A (\pi - \pi_0) + (\bX_j - \pi)^\top A (\pi - \pi_0) + ||A^{1/2}(\pi - \pi_0)||^2_2 \Big\}.
\end{align*}
Accordingly, we have $\mE \left[ U_{\text{quad}}\right] = 0$ and $\mE \left[ U_{\text{linear}} \right] = || A^{1/2}(\pi - \pi_0)||_2^2$. Under (S.1), $U_A$ is dominated by $U_{\text{linear}}$ and thus
\begin{align*}
\frac{U_A - || A^{1/2}(\pi - \pi_0) ||_2^2  }{\sqrt{\text{Var}(U_A) }} = \frac{U_{\text{linear}} - || A^{1/2}(\pi - \pi_0) ||_2^2  }{ \sqrt{\text{Var} (U_{\text{linear}})}} + o_P(1),
\end{align*}
whereas under (S.2), $U_A$ is dominated by $U_{\text{quad}}$ so that
\begin{align*}
\frac{U_A - ||A^{1/2}(\pi - \pi_0)||_2^2  }{\sqrt{\text{Var}(U_A) }} = \frac{U_{\text{quad}}}{ \sqrt{\text{Var} (U_{\text{quad}} )}} + o_P(1).
\end{align*}
Hence, in order to establish the asymptotic normality of $U_A$, we need to study the limiting behavior of $U_{\text{linear}}$ and $U_{\text{quad}}$ under each scenario. The result is summarized in the following theorem.

\begin{theorem}[Asymptotic normality of $U_A$ under the alternative] \label{Thm: Gaussian Limit under alternative}
Assume either i) \emph{(S.1)} and $(\pi-\pi_0)^\top A \Sigma A(\pi - \pi_0) < \infty$, or ii) \emph{(S.2)} and the condition (\ref{Eq: Gaussian Martingale Assumptions}) given in Theorem \ref{Thm: Global Gaussian}.
Then
\begin{align*}
\frac{U_A - || A^{1/2}(\pi - \pi_0) ||_2^2  }{\sqrt{\text{\emph{Var}}(U_A) }} \convD \mathcal{N}(0,1),
\end{align*} 
where 
\begin{align*}
\text{\emph{Var}}\left( U_A \right) = \binom{n}{2}^{-1} \Big\{ \emph{tr} \{ (A \Sigma)^2 \}+ 2(n-1) (\pi-\pi_0)^\top A \Sigma A(\pi - \pi_0) \Big\}. 
\end{align*}
\end{theorem}

Theorem \ref{Thm: Gaussian Limit under alternative} together with Theorem \ref{Thm: Global Gaussian} allows us to describe the power function of $U_A$ under the Gaussian asymptotic. Let $z_\alpha$ be the upper $\alpha$ quantile of the standard normal distribution. For notational simplicity, let us denote 
\begin{align*}
\Lambda_0 = \text{tr} \{ (A \Sigma_0)^2 \}, \ \ \Lambda_1 = \text{tr} \{ (A \Sigma)^2 \} \quad \text{and} \quad \Lambda_2 = 2(n-1) (\pi-\pi_0)^\top A \Sigma A(\pi - \pi_0),
\end{align*}
where $\Sigma_0 $ and $\Sigma$ are the covariance matrix of $\bX$ under the null and the alternative hypothesis, respectively. Then the power is approximated by
\begin{align} \label{Eq: power approximation via asymptotic normality}
\beta_{n,d}(\pi_0, \pi_1, A) = \Phi \left(  -\frac{\sqrt{\Lambda_0}}{\sqrt{\Lambda_1 + \Lambda_2}}z_{\alpha}  + \sqrt{\binom{n}{2}} \frac{||A^{1/2}(\pi-\pi_0)||_2^2}{\sqrt{\Lambda_1 + \Lambda_2}}  \right) + o(1).
\end{align}
Under (S.1) together with the additional assumption (S.3) below: 
\begin{enumerate}\setlength{\itemindent}{.2in}
\item [\textbf{(S.3)}] \ $n^{-1}\text{tr}( (A\Sigma_0)^2)=  o\big((\pi-\pi_0)^\top A \Sigma A(\pi - \pi_0) \big)$,
\end{enumerate}
the power function of $U_A$ can be further simplified to
\begin{align*}
\beta_{n,d}(\pi_0, \pi_1, A)  = \Phi \left(  \frac{\sqrt{n}||A^{1/2}(\pi-\pi_0)||_2^2}{\sqrt{ 4 (\pi-\pi_0)^\top A \Sigma A(\pi - \pi_0) }}  \right) + o(1).
\end{align*}
On the other hand, under (S.2), the approximation becomes
\begin{align*}
\beta_{n,d}(\pi_0, \pi_1, A)  = \Phi \left(   -\frac{\sqrt{\text{tr} \{ (A \Sigma_0)^2 \}}}{\sqrt{\text{tr} \{ (A \Sigma)^2 \}}}z_\alpha + \frac{n ||A^{1/2}(\pi-\pi_0)||_2^2}{\sqrt{2 \text{tr} \{ (A \Sigma)^2 \} }}  \right) + o(1).
\end{align*}

\vskip 2em

\section{Minimax Optimality} \label{Sec: Minimax Upper Bound}
As discussed before, $\chi^2$ statistic tends to have a large variance by putting too much weight on small entries of $\pi_0$. Consequently, the resulting test can perform poorly and is not minimax optimal in the high-dimensional setting \citep{balakrishnan2018hypothesis}. This motivates us to consider different weights for the test statistic. In this section, we discuss the choice of the weight matrix $A$ from a minimax point of view. To formulate the minimax problem, we modify the hypotheses given in (\ref{Eq: Hypothesis}) as
\begin{align} \label{Eq: L1 Hypothesis}
H_0:  \pi = \pi_0 \quad \text{versus} \quad H_1: ||\pi - \pi_0 ||_1 \geq \epsilon_n,
\end{align}
where $||x ||_1$ is the $\ell_1$ norm of $x \in \mathbb{R}^d$. Let us consider a set of level $\alpha$ test functions, $\phi : \{\bX_i\}_{i=1}^n \mapsto \{0,1\}$, such that
\begin{align}
\Phi_{n,\alpha} = \Big\{ \phi : \mP_{H_0}^n \left(\phi = 1 \right) \leq \alpha, 0 < \alpha < 1 \Big\}.
\end{align}
Then the global minimax risk \citep[see e.g.,][]{valiant2017automatic,balakrishnan2017hypothesis} is defined as the supremum over the local minimax risk:
\begin{align*}
R_n(\epsilon_n) = \sup_{\pi_0 \in \Omega} R_n(\epsilon_n,  \pi_0),
\end{align*}
where the local minimax risk is given by
\begin{align*}
R_n(\epsilon_n, \pi_0) = \inf_{\phi \in \Phi_{n,\alpha}} \sup \Big\{ \mE_{H_1}[1- \phi] : ||\pi- \pi_0||_1 \geq \epsilon_n, \pi \in \Omega  \Big\}.
\end{align*}
For a given $\delta \in (0, 1-\alpha)$, the global minimum separation rate is characterized by
\begin{align*}
\epsilon_n^\ast = \inf \Big\{ \epsilon_n :  R_n(\epsilon_n) \leq \delta \Big\}.
\end{align*}
Under the given setting, \cite{valiant2017automatic} show that the global minimax rate is
\begin{align*}
\epsilon_n^\ast \asymp  \frac{d^{1/4}}{\sqrt{n}}.
\end{align*}

The main objective of this section is to find a sufficient condition for $A$ that results in minimax rate optimal test based on the $U$-statistic. We first describe that the test based on the $U$-statistic with a mixture weight is minimax rate optimal in Section~\ref{Section: $U$-statistic weighted by a mixture distribution} and we go on to generalize this result in Section~\ref{Section: Generalization}.

\vskip .8em

\subsection{$U$-statistic weighted by a mixture distribution} \label{Section: $U$-statistic weighted by a mixture distribution}
The weight used in $U_{\pi_0}$ often results in a high variance of the test statistic especially when $\pi_0$ is sparse. $U_I$ does not suffer from the same variance issue but its weight does not use information of the null distribution. We combine $U_{\pi_0}$ and $U_{I}$ to reduce the disadvantages associated with each and obtain a minimax rate optimal test. Let us define the mixture distribution by
\begin{align}\label{Eq: Mixture distribution}
\pi_{\text{mix}} = \frac{1}{2} \pi_0  + \frac{1}{2} \pi_{\text{unif}},
\end{align}
where $\pi_{\text{unif}} = (1/d,\ldots,1/d)$. Then by using $A_{\text{mix}} = \text{diag} \{ \pi_{\text{mix},1}^{-1}, \ldots, \pi_{\text{mix},d}^{-1}  \}$, the resulting $U$-statistic is defined by
\begin{align}  \label{Eq: U-statistic with mixture weight}
U_{\text{mix}} = \binom{n}{2}^{-1}  \sum_{1 \leq i < j \leq n} \left(\bX_i - \pi_0 \right)^\top A_{\text{mix}} (\bX_j - \pi_0).
\end{align}
To test (\ref{Eq: L1 Hypothesis}), we reject the null hypothesis when $U_{\text{mix}}$ is greater than the critical value:
\begin{align*}
\phi(U_{\text{mix}} ) & ~=~ I \left( U_{\text{mix}} >  \sqrt{ \frac{1}{\alpha}\binom{n}{2}^{-1}  \text{tr} \{ (A_{\text{mix}} \Sigma_0)^2 \}}  \right),
\end{align*}
where $\Sigma_0 = \text{diag}\left(\pi_0 \right) - \pi_0 \pi_0^\top$. Then we can see that $\phi(U_{\text{mix}} )$ has nontrivial power when $\epsilon_n \asymp d^{1/4} / \sqrt{n}$ and thus it is global minimax rate optimal. We formally state this result in Theorem~\ref{Thm: Global minimax optimality of $U_{w}$} which holds for more general test statistics. 

%\begin{theorem}[Global minimax optimality of $U_{\text{mix}}$] \label{Thm: Global Minimax Optimality based on a mixture distribution}
%	For testing (\ref{Eq: L1 Hypothesis}), the test based on $\phi(U_{\text{mix}})$ has size at most $\alpha$. In addition, suppose there exists a universal constant $C>0$ independent of $n$ and $d$ such that 
%	\begin{align}
%	\epsilon_n^2 \geq \frac{C\sqrt{d}}{n} \left[ \frac{1}{\sqrt{\alpha}} + \frac{1}{\zeta} \right],
%	\end{align}
%	for any $\zeta \in (0,1]$, then we have $\mP_{H_1}(\phi(U_{\text{mix}}) = 0) \leq \zeta$. Hence, $\phi(U_{\text{mix}})$ is global minimax optimal.
%\end{theorem}

\vskip .8em

\begin{remark} \label{Remark: Collision-based test}
	\cite{diakonikolas2016collision} show that the collision-based test statistic 
	\begin{align*} 
	W = \sum_{1 \leq i < j \leq n} \emph{\bX}_i^\top \emph{\bX}_j
	\end{align*}
	is minimax rate optimal for multinomial uniformity testing. For the uniform null case, $U_{\text{mix}}$ is equivalent to the collision-based test statistic $W$; thereby, our result can be viewed as a generalization of \cite{diakonikolas2016collision} to arbitrary null probabilities. 
\end{remark}

\begin{remark}
	Poissonization, where the sample size has a Poisson distribution, is a standard assumption in the literature to construct the upper bound of the minimax risk. Under Poissonization, several statistics have been proposed to obtain the minimax optimality \citep{valiant2017automatic,balakrishnan2017hypothesis}. We would like to emphasize that our minimax result is established without assuming Poissonization.
\end{remark}

\vskip 1em

\subsection{Generalization} \label{Section: Generalization}
The mixture distribution in (\ref{Eq: Mixture distribution}) can be generalized by considering an arbitrary but fixed $\gamma \in (0,1)$ such that
\begin{align*}
\pi_{\text{mix}}^{(\gamma)} = \gamma \pi_0   + (1 - \gamma) \pi_{\text{unif}}.
\end{align*}
For a given weight vector $w \in \mathbb{R}^d$, we say that $w \in \mathbb{R}^d$ is \emph{comparable} to $\pi_{\text{mix}}^{(\gamma)}$, if there exist fixed constants $C_1, C_2  >0$ independent of $n$ and $d$ such that $C_1 \pi_{\text{mix},i}^{(\gamma)} \leq w_i \leq C_2 \pi_{\text{mix},i}^{(\gamma)}$ for all $i=1,\ldots, d$. We denote a weight vector comparable to $\pi_{\text{mix}}^{(\gamma)}$ by $$w \sim \pi_{\text{mix}}^{(\gamma)}.$$ Based on these notations, let us define a class of weight matrices:
\begin{align}
\mathcal{A}_w = \Big\{ \text{diag}\left( w^{-1} \right) \in \mathbb{R}^{d \times d } : w \in \mathbb{R}^d ~ \text{and}  ~ w \sim \pi_{\text{mix}}^{(\gamma)} \ \  \text{for some} \ \gamma \in (0,1) \Big\}.
\end{align}
Then the test based on the $U$-statistic associated with any $A_w \in \mathcal{A}_w$: 
\begin{align*}
\phi(U_{w} ) & ~=~ I \left( U_{w} >  \sqrt{ \frac{1}{\alpha}\binom{n}{2}^{-1}  \text{tr} \{ (A_{w} \Sigma_0)^2 \}}  \right)
\end{align*}
is global minimax rate optimal. The result is summarized in the next theorem.

\begin{theorem}[Global minimax optimality of $U_{w}$] \label{Thm: Global minimax optimality of $U_{w}$}
	For testing (\ref{Eq: L1 Hypothesis}), the test based on $\phi(U_{w})$ has size at most $\alpha$. In addition, suppose there exists a universal constant $C>0$ independent of $n$ and $d$ such that 
	\begin{align} \label{Eq: Condition for critical radius}
	\epsilon_n^2 \geq \frac{C\sqrt{d}}{n} \left[ \frac{1}{\sqrt{\alpha}} + \frac{1}{\zeta} \right],
	\end{align}
	for any $\zeta \in (0,1]$, then we have $\mP_{H_1}(\phi(U_{w}) = 0) \leq \zeta$. Hence, $\phi(U_{w})$ is global minimax optimal.
\end{theorem}

\vskip .8em

Here we provide several examples that belong to the proposed framework. 

\begin{example}[Truncated $\chi^2$] \label{Example: Truncated weights}
	\cite{balakrishnan2017hypothesis} show that the test based on the truncated $\chi^2$ test statistic is global minimax rate optimal. Unlike the classical $\chi^2$ statistic, the truncated $\chi^2$ test statistic is weighted by $\theta_{\text{\emph{trunc}},j} = \max \{ \pi_{0,j}, 1/d \}$ for $j = 1,\ldots, d$. Note that $\theta_{\text{\emph{trunc}}} \sim \pi_{\text{\emph{mix}}}$ since $\pi_{\text{\emph{mix}},j} \leq \theta_{\text{\emph{trunc}},j} \leq 2 \pi_{\text{\emph{mix}},j}$
	for all $j = 1,\ldots, d$. Therefore, it satisfies the comparable condition with $C_1 = 1$ and $C_2 = 2$.
\end{example}

\begin{example}[$\ell_p$-type mixture]
	For $p \geq 1$, let us define
	\begin{align*}
	\theta_{\ell_p,j} = \left( \frac{\pi_{0,j}^p + \pi_{\text{\emph{unif}},j}^{p} }{2}  \right)^{1/p}
	\end{align*}
	for $j=1,\ldots,d$. Then we observe that $\theta_{\ell_p} \sim \pi_{\text{\emph{mix}}}$ since $\pi_{\text{\emph{mix}},j} \leq \theta_{\ell_p,j} \leq 2^{1 - 1/p} \pi_{\text{\emph{mix}},j}$ for all $j=1,\ldots, d$, where we used $||x||_1 \leq 2^{1 - 1/p} ||x||_p  \leq 2^{1 - 1/p} ||x||_1$ for $p \geq 1$. In fact, if $p=\infty$, it corresponds to the truncated weight as $\theta_{\ell_\infty} =\theta_{\text{\emph{trunc}}}$.
\end{example}

\vskip .8em

\section{Simulations} \label{Sec: Simulations}
In this section, we provide numerical results to illustrate the finite sample performance of the proposed methods. In the first simulation study, we compare power between Pearson's chi-squared test and the proposed tests based on $U$-statistics. We consider four different $U$-statistics: $U_{\pi_0}$, $U_{I}$, $U_{mix}$ and $U_{trunc}$ where $U_{trunc}$ is the $U$-statistic with truncation weights described in Example~\ref{Example: Truncated weights}. We let the null distribution $\pi_0$ have a power law distribution where the probability of the $i$th bin is proportional to the $r_0$th power of its index, i.e. $\pi_{0,i} \propto i^{r_0}$ for $i=1,\ldots,d$. When $r_0$ is close to zero, then the null distribution becomes close to the uniform distribution. On the other hand, when $r_0$ has a large value, the null distribution becomes skewed to the left.  %When the null distribution is uniformly distributed, i.e. $r_0 = 0$, the considered test statistics become equivalent to each other (up to some scaling and location factors), resulting in the same power performance. For this reason, 
In our simulations, we consider two null distributions with $r_0 = 1$ and $r_0=5$. The alternative distribution $\pi$ is also chosen to have a power law distribution as $\pi_i \propto i^{r}$ for $i=1,\ldots,d$ and we change $r$ to describe different power behaviors. The null and alternative distribution of each statistic are estimated via Monte Carlo simulations with 1000 repetitions where we take the sample size and the number of bins as $n=200$ and $d=2000$, respectively.

\begin{figure}[!t]
	\subfigure{\includegraphics[width=7.7cm]{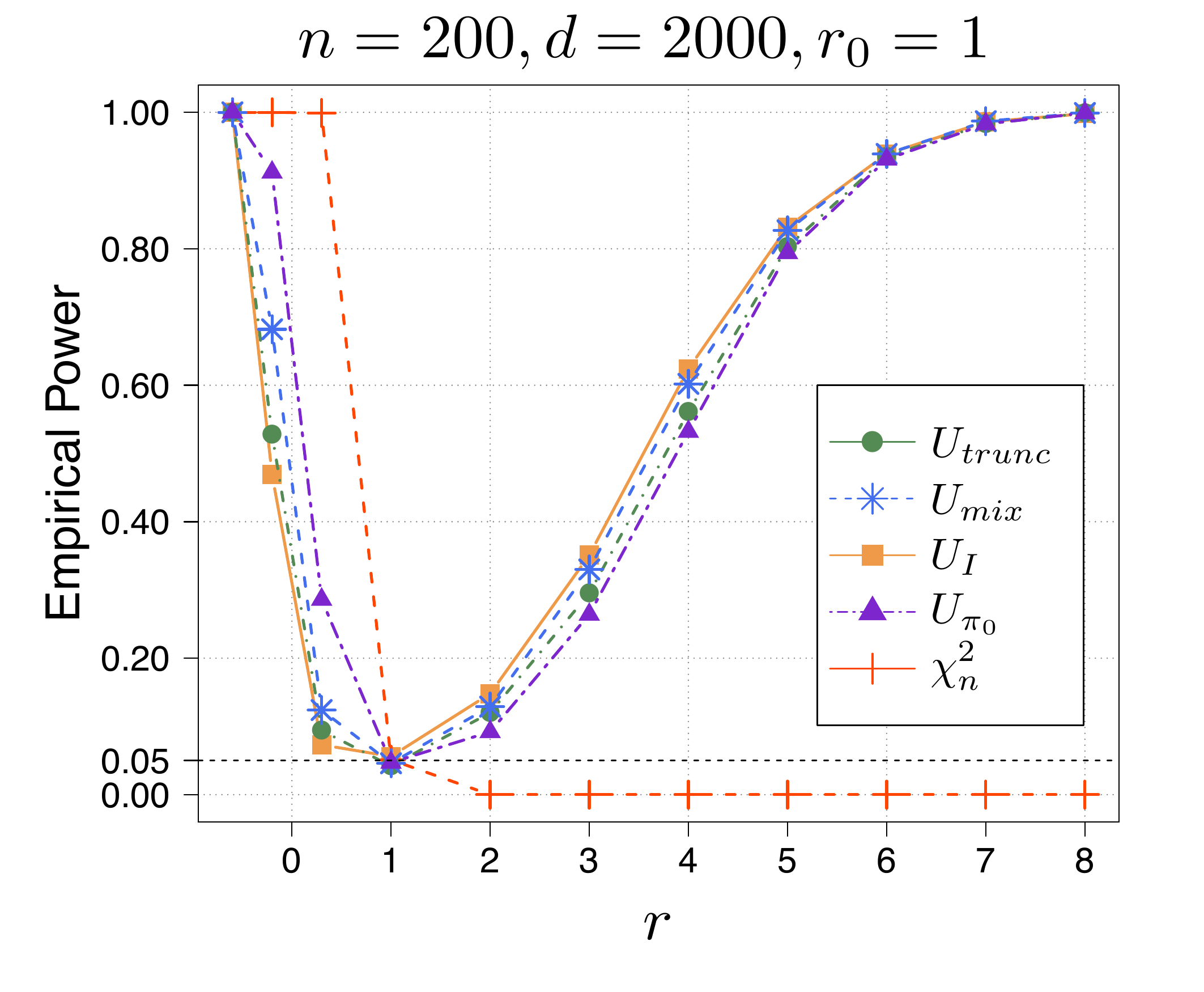}}
	\hfill
	\subfigure{\includegraphics[width=7.7cm]{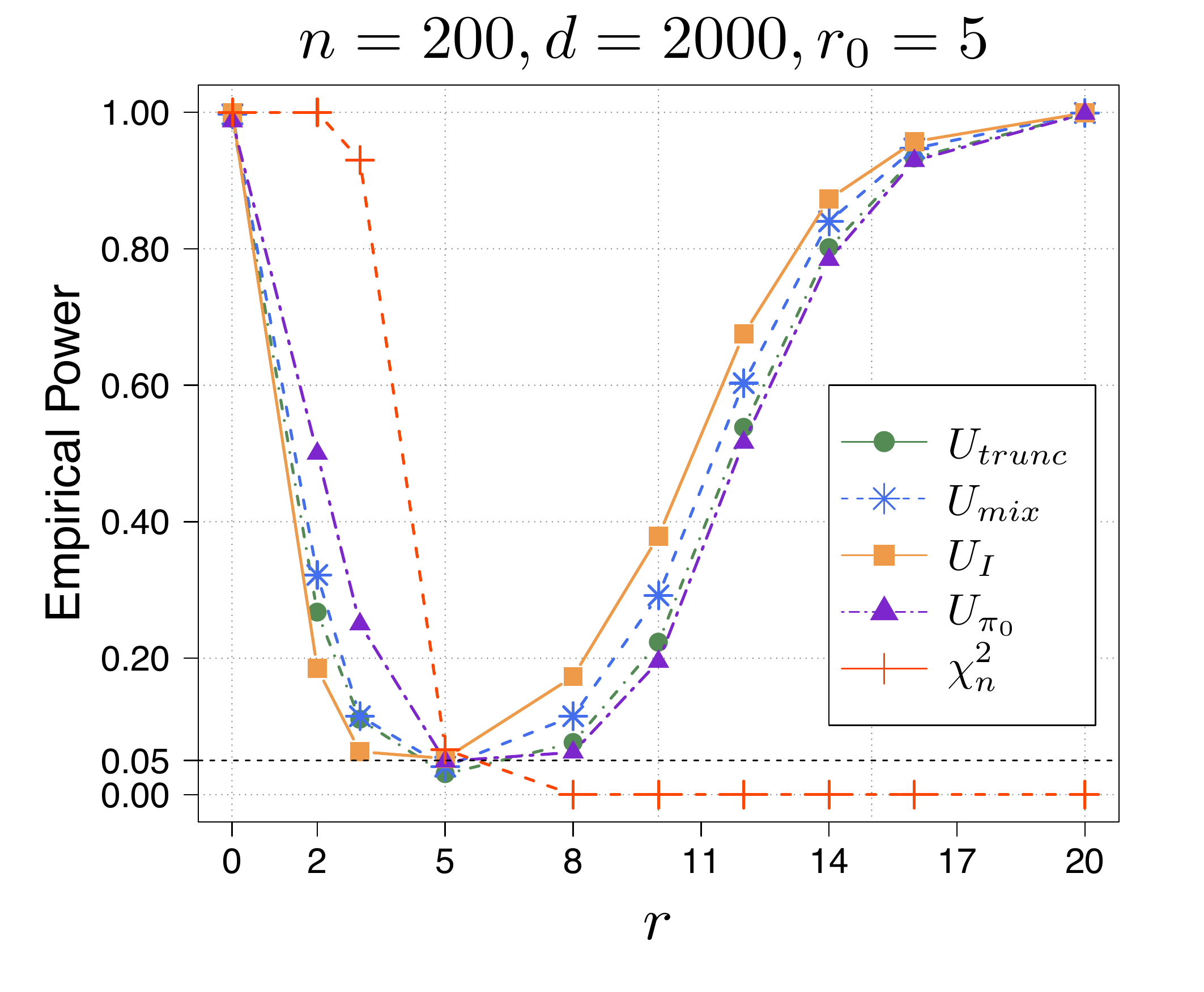}}
	\hfill
	\caption{Power comparisons between five different tests based on $U_{trunc}$, $U_{mix}$, $U_I$, $U_{\pi_0}$ and Pearson's $\chi_n^2$ at significance level $\alpha=0.05$.}  \label{Figure: Power comparison}
\end{figure}

The simulation results are presented in Figure~\ref{Figure: Power comparison}. From the results, we observe that Pearson's $\chi_n^2$ test shows entirely different behaviors between two alternatives where $(i) ~ r_0 < r$ and $(ii)~r_0 > r$. Specifically, when $r_0<r$, Pearson's $\chi_n^2$ test is extremely biased and has zero power to reject the null hypothesis, whereas it has the highest power among the considered tests when $r_0 > r$. In contrast, the tests based on the $U$-statistics are considerably robust toward the testing bias and perform reasonably well against the entire range of alternatives. This illustrates the benefit of the proposed $U$-statistic framework under the  high-dimensional regime. For the comparison between the $U$-statistics, no test is uniformly more powerful than the others. In particular, the test based on $U_{\pi_0}$ outperforms the other tests when $r_0 > r$, but underperforms when $r_0 < r$. On the other hand, the test based on $U_I$ performs the best when $r_0 < r$ and performs the worst when $r_0 > r$. The test based on $U_{mix}$ are usually the second best and better than the one based on $U_{trunc}$.

\begin{figure}[!t]
	\subfigure{\includegraphics[width=7.7cm]{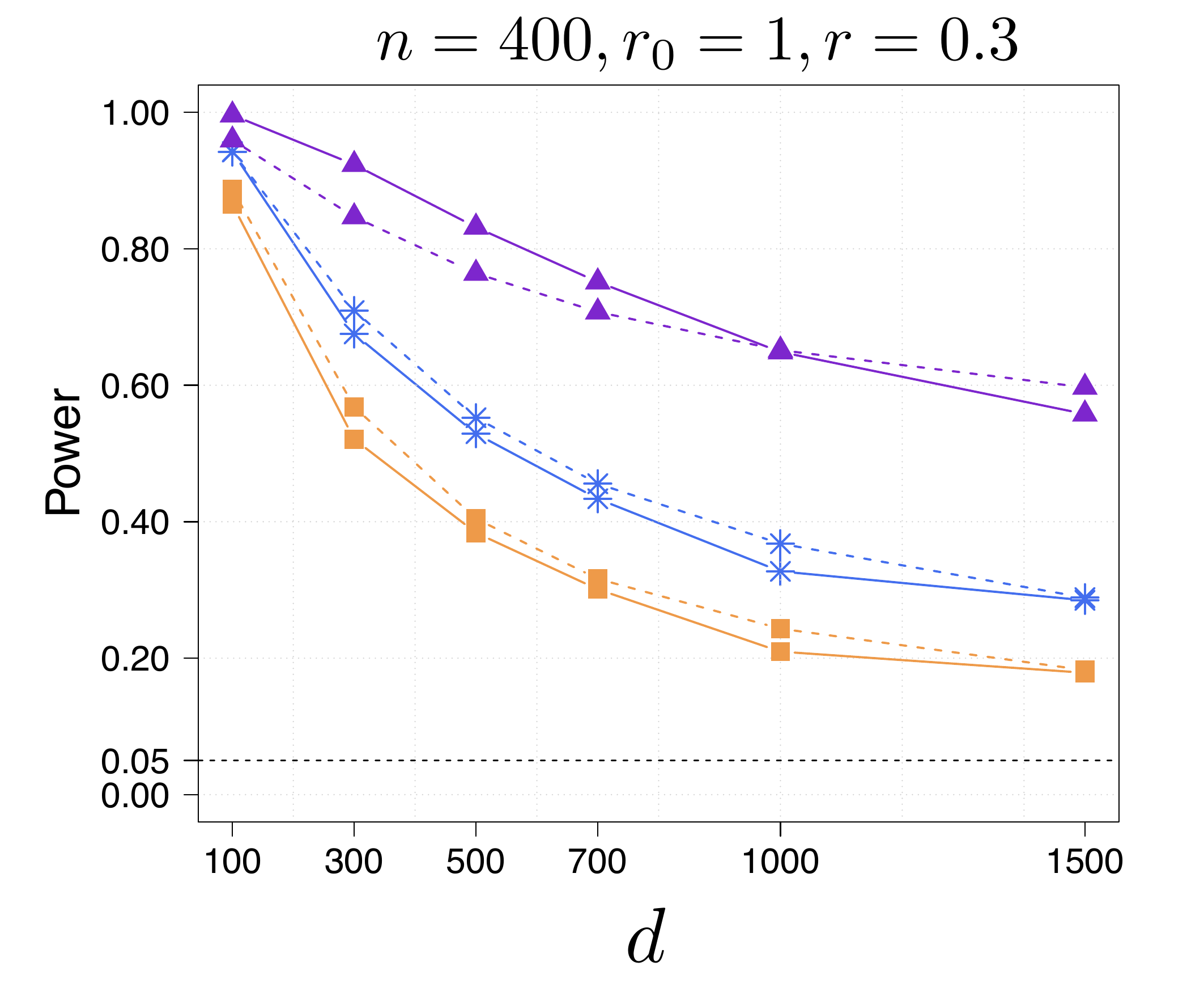}}
	\hfill
	\subfigure{\includegraphics[width=7.7cm]{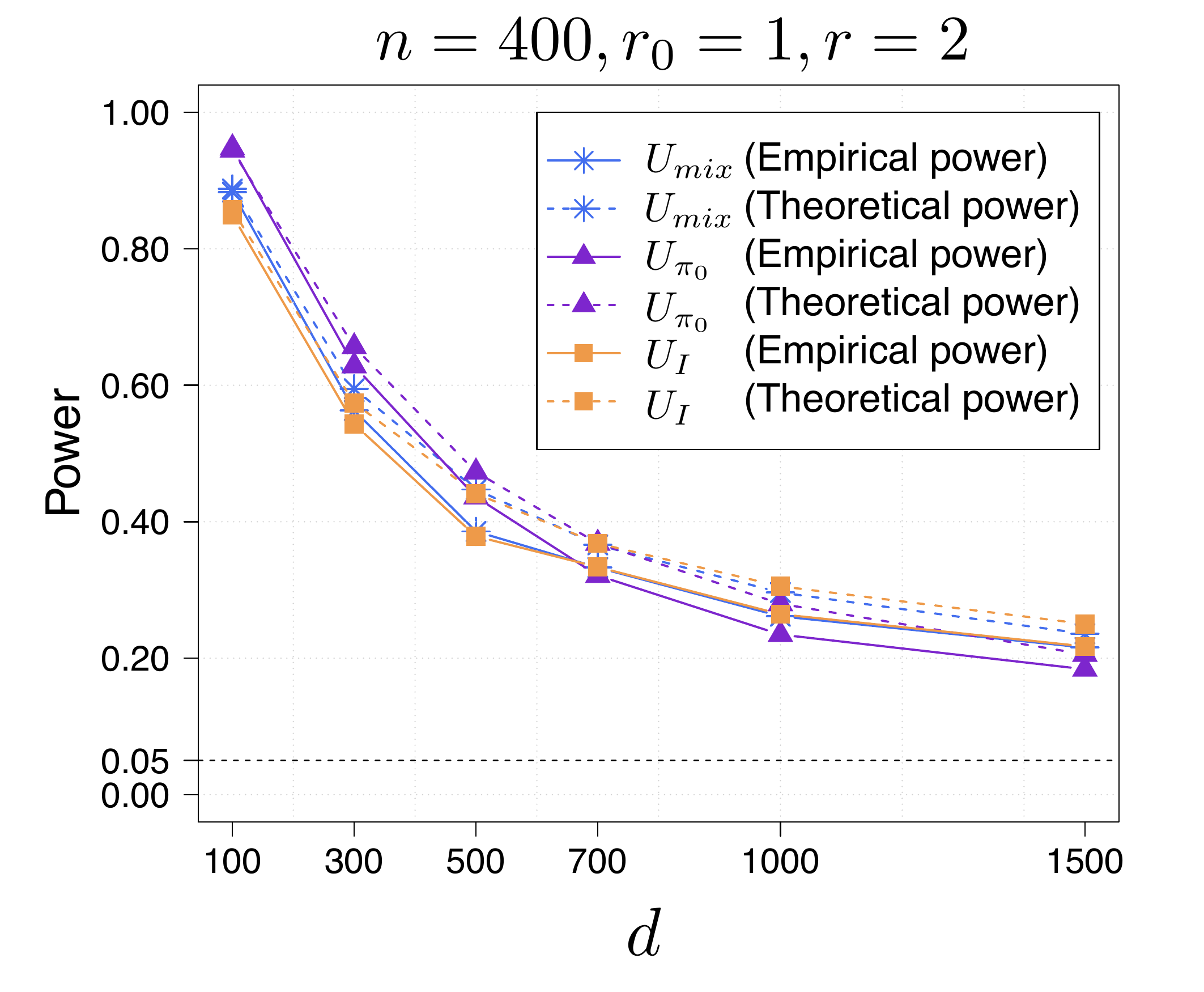}}
	\hfill
	\caption{Comparisons between the empirical power and the theoretical power based on the normal approximation at significance level $\alpha =0.05$.} \label{Figure: Power approximation}
\end{figure}

In the second simulation study, we compare the empirical power of the tests and the corresponding theoretical power based on the normal approximation established in (\ref{Eq: power approximation via asymptotic normality}). For the comparison, we consider three $U$-statistics: $U_{\pi_0}$, $U_{I}$ and $U_{mix}$, and choose the null distribution as $\pi_{0,i} \propto i$ for $i=1,\ldots,d$. Under the alternative, we consider two power law distributions with $r=0.3$ and $r=2$ where $\pi_i \propto i^r$ for $i=1,\ldots,d$. As before, the null and alternative distribution of each statistic are estimated via Monte Carlo simulations with 1000 repetitions where we take the sample size and the number of bins as $n=400$ and $d=100,300,500,700,1000,1500$. The results are given in Figure~\ref{Figure: Power approximation}. It is seen from the results that the power approximation via asymptotic normality looks fairly robust over different dimensions especially against the alternative distribution with $r=2$. %Although the power approximation may not be accurate in other scenarios, we believe that it would be still useful for large $n$ and $d$ by providing a rough idea of statistical power of the test without relying on heavy simulations.

\vskip .8em

\section{Summary and Discussion} \label{Sec: Conclusions}
In this work, we introduced a family of $U$-statistics for multinomial goodness-of-fit tests and investigated their asymptotic behaviors in the high-dimensional regime. Specifically, we established the conditions under which the $U$-statistic is approximately Poisson or Gaussian, and studied its power function under each asymptotic regime. We also proposed a class of weights for the $U$-statistic and showed the minimax optimality of the resulting tests. Despite the fact that the proposed tests achieve minimax rate optimality, they still have room for improvement. In particular, the considered class of weight functions only uses the information of $\pi_0$ but not $\pi$. When prior information about $\pi$ is available (e.g. differences exist in specific bins with high probability), then it is possible to design more powerful test by incorporating that information. In this case, it would be beneficial to use our asymptotic results and choose $A$ that maximizes the asymptotic power function under the given restrictions. We reserve this topic for future work.

\section*{Acknowledgements} The author would like to thank Sivaraman Balakrishnan and Larry Wasserman for their valuable comments which helped to improve the manuscript.

\bibliographystyle{apalike}
\bibliography{reference}

%\clearpage

\appendix
\allowdisplaybreaks

\section{Proofs} \label{Sec: Appendix}
\subsection{Proof of Lemma \ref{Lemma: Pearson V-stat}}
We will provide a more general result by considering an arbitrary positive diagonal matrix $A = \text{diag}\{a_1, \ldots, a_d\}$ in the kernel. In other words, we will show the following holds:
\begin{align} \label{Eq: V-stat1}
\sum_{j=1}^d a_j\left( Y_j - n\pi_{0,j} \right)^2 = \sum_{1 \leq i \leq n} \sum_{1 \leq j \leq n} (\bX_i - \pi_0)^\top A (\bX_j - \pi_0).
\end{align}
Then the result follows by setting $a_i = \left( n\pi_{0,i}\right)^{-1}$. 

\vskip .5em

First, we decompose the left hand side of (\ref{Eq: V-stat1}) into the three parts: 
\begin{align*}
\sum_{j=1}^d a_j\left( Y_j - n\pi_{0,j} \right)^2 = \underbrace{ \sum_{j=1}^d a_j Y_j^2}_{(i)} - \underbrace{2 n\sum_{j=1}^d  a_j Y_j \pi_{0,j}}_{(ii)} + \underbrace{n^2\sum_{j=1}^d a_j \pi_{0,j}^2}_{(iii)} 
\end{align*}
and treat them separately. 

\vskip .5em

\noindent \textbf{Part (i).} Recall that $Y_j = \sum_{i=1}^n I(X_{i,j}=1)$, and thus
\begin{align*}
\sum_{j=1}^d a_j  \left( \sum_{i=1}^n I(X_{i,j}=1) \right)^2 & =  \sum_{j=1}^d a_j \left[ \sum_{i=1}^n I(X_{i,j}=1) + 2 \sum_{i < i^\prime} I(X_{i,j}=1) I(X_{i^\prime, j}=1)  \right] \\
& = \sum_{i=1}^n \sum_{j=1}^d a_j I(X_{i,j}=1)  + 2 \sum_{i < i^\prime} \sum_{j=1}^d a_j I(X_{i,j}=1) I(X_{i^\prime, j}=1) \\
& = \sum_{i=1}^n \bX_i^\top A \bX_i + 2 \sum_{i < i^\prime} \bX_i^\top A \bX_{i^\prime}   = \sum_{i=1}^n \sum_{i^\prime=1}^n \bX_i^\top A \bX_{i^\prime}.
\end{align*}

\noindent \textbf{Part (ii).} Similar to the first part, 
\begin{align*}
2 n\sum_{j=1}^d  a_j Y_j \pi_{0,j} & = 2n \sum_{j=1}^d a_j \left(   \sum_{i=1}^n I(X_{i,j}=1) \right) \pi_{0,j} \\ 
& = 2n \sum_{i=1}^n \sum_{j=1}^d a_j I(X_{i,j}=1) \pi_{0,j} = 2n \sum_{i=1}^n \bX_i^\top A \pi_0 \\
& = \sum_{i=1}^n \sum_{i^\prime=1}^n \left( \bX_i^\top A \pi_0 + \bX_{i^\prime}^\top A \pi_0 \right).
\end{align*}
\noindent \textbf{Part (iii).} The last part is straightforward, 
\begin{align*}
n^2 \sum_{j=1}^d a_j \pi_{0,j}^2 =  \sum_{i=1}^n \sum_{i^\prime=1}^n \pi_0^\top A \pi_0.
\end{align*} 
Combining the three parts, we can get the desired result.

\vskip .8em

\subsection{Proof of Theorem~\ref{Thm: Poisson approximation}}

The proof is mainly based on Theorem 1 of \cite{arratia1989two}. We describe it in Theorem~\ref{Thm: Arratia theorem 1}. Before we get to the main proof, we provide several lemmas. 

\vskip .8em 

Recall the following decomposition of $U_I$:
\begin{align*}
U_I = \binom{n}{2}^{-1} \underbrace{\sum_{1 \leq i < j \leq n} \bX_i^\top \bX_j}_{W} - \frac{2}{n} \sum_{i=1}^n \left( \bX_i^\top \pi_0 + \pi_0^\top \pi_0\right).
\end{align*}

\vskip .8em

Suppose $\binom{n}{2} \pi_0^\top \pi_0 \rightarrow \eta$ as $n,d \rightarrow \infty$. As preliminary results, we are interested in the conditions that result in
\begin{align*}
	\quad W= \sum_{1 \leq i < j \leq n} \bX_i^\top \bX_j  \convD \text{Pois}(\eta) \quad \text{where }   \eta \in (0, \infty).
\end{align*}

\noindent Note that $W$ is the sum of locally dependent indicator variables. The limiting distribution of $W$ has been studied under the name of the birthday problem \citep[see e.g.,][]{dasgupta2005matching}. Let $Z$ be a Poisson random variable with $\mE[W] = \mE[Z]$. Here, our interest is in the total variation distance between $W$ and $Z$, i.e.
\begin{align*}
d_{TV} (W,Z) = 2 \sup_{A \in \mathbb{Z}^{+}} \Big| \mP (W \in A) - \mP (Z \in A ) \Big|.
\end{align*}
In order to bound the total variation distance, we employ Theorem 1 of \cite{arratia1989two}:

\begin{theorem}[Theorem 1 of \cite{arratia1989two}] \label{Thm: Arratia theorem 1}
	Let $\mathcal{I}$ be an arbitrary index set, and for $\alpha \in \mathcal{I}$, let $X_\alpha$ be a Bernoulli random variable with $p_\alpha = \mP(X_\alpha= 1)>0$. Let $K = \sum_{\alpha \in \mathcal{I}} X_\alpha$ and $\mE[K] = \sum_{\alpha \in \mathcal{I}} p_\alpha \in (0, \infty)$. For each $\alpha \in \mathcal{I}$, suppose we have chosen $B_\alpha \in \mathcal{I}$ with $\alpha \in B_\alpha$. Define 
	\begin{align*}
	& b_1 = \sum_{\alpha \in \mathcal{I}} \sum_{\beta \in B_\alpha} p_\alpha p_\beta, \ \  b_2 = \sum_{\alpha \in \mathcal{I}} \sum_{\alpha \neq \beta \in B_\alpha} \mE[X_\alpha X_\beta],  \ \ \text{and}  \ \  b_3 = \sum_{\alpha \in \mathcal{I}} \mE \bigg|  \mE \bigg\{ X_\alpha - p_\alpha \bigg|  \sum_{\beta \in \mathcal{I} - B_\alpha } X_\beta  \Big\}  \bigg|.
	\end{align*} 
	Let $Z$ be a Poisson random variable with $\mE[K] = \mE[Z] = \eta$. Then
	\begin{align} \label{Eq: Arriatia main result}
	d_{TV}(K,Z) \leq 2 \left[ (b_1 + b_2) \frac{1-e^{-\eta}}{\eta} + b_3  (1 \wedge 1.4 \eta^{-1/2})\right].
	\end{align}
\end{theorem}

\vskip .8em

As a corollary of Theorem \ref{Thm: Arratia theorem 1}, we present the total variation distance between $W$ and the Poisson random variable $Z$.

\begin{corollary} \label{Corollary:PoissonLimit}
	Let $Z$ be a Poisson random variable with $\mE[Z] = \mE[W] = \binom{n}{2} \pi^\top \pi= \eta_{n,d}$. Then 
		\begin{align*}
	d_{TV}(W, Z) & ~ \leq  ~ 2n^3\frac{1 - e^{-\eta_{n,d}}}{\eta_{n,d}} \left[ \sum_{i=1}^d \pi_i^3 + \Bigg\{ \sum_{i=1}^d \pi_i^2\Bigg\}^2  \right]. 
	\end{align*}
	\begin{proof}
		Denote $p_1 = \mP\left(\bX_1^\top \bX_2=1 \right)$ and $p_{2} = \mP \left(\bX_1^\top \bX_2 \bX_1^\top \bX_3=1\right)$. In addition, let $\mathcal{I}$ be a set of all indices $(i,j)$ where $1 \leq i < j \leq n$ so that $|\mathcal{I}| = \binom{n}{2}$ and $\mathcal{B}_i$ be a collection of $\bX_j^\top \bX_k$ that is dependent on $\bX_i$. Here $|\mathcal{B}_i | = \binom{n}{2} - \binom{n-2}{2}= 2n -3$ for all $i$.
		
		\vskip .5em
		
		\noindent Thus we have 
		\begin{align*}
		b_1 & = p_1^2 |\mathcal{I} | | \mathcal{B}_i| = \Bigg\{ \sum_{i=1}^d \pi_i^2 \Bigg\}^2 \binom{n}{2}(2n-3) \leq \Bigg\{ \sum_{i=1}^d \pi_i^2 \Bigg\}^2 n^3,  \\[.5em]
		b_2 & = p_{2} |\mathcal{I} | \left( | \mathcal{B}_i | - 1 \right) = \sum_{i=1}^d \pi_i^3 \binom{n}{2} (2n-4) = \sum_{i=1}^d \pi_i^3 n(n-1)(n-2)  \leq  \sum_{i=1}^d \pi_i^3 n^3.
		\end{align*}
		and $b_3 = 0$ by the construction of $\mathcal{B}_i$. Then the proof is complete by applying Theorem \ref{Thm: Arratia theorem 1}.
	\end{proof}
\end{corollary}

So far, we have investigated the asymptotic results of $W$. Now, let us turn our attention to $U_I$. Recall that $U_I$ and $W$ are related as
\begin{align*}
\binom{n}{2}U_I =  W - (n-1)\sum_{i=1}^n \bX_i^\top \pi_0 + \binom{n}{2} \pi_0^\top \pi_0.
\end{align*}
Hence, in order to attain the Poisson approximation for $\binom{n}{2} U_I$, one might need to control the last two terms properly.
Note that
\begin{align*}
\mE \left[  (n-1)\sum_{i=1}^n \bX_i^\top \pi_0 \right] & = 2 \binom{n}{2} \pi^\top \pi_0 = 2\eta_2 + o(1), \\
\text{Var} \left[ (n-1)\sum_{i=1}^n \bX_i^\top \pi_0\right]  & = (n-1)^2 n \left[ \mE\left[ \bX_1^\top \pi_0 \bX_1^\top \pi_0 \right]  - \left( \mE\left[\bX_1^\top \pi_0\right] \right)^2 \right] \\
& \leq n^3 \left( \sum_{i=1}^d \pi_i \pi_{0,i}^2 -  \left(\sum_{i=1}^d \pi_{i} \pi_{0,i} \right)^2  \right).
\end{align*}
Hence, under (P.2) and (P.3), Chebyshev's inequality yields
\begin{align} \label{Eq: Aim1 - Poisson limit}
(n-1)\sum_{i=1}^n \bX_i^\top \pi_0 \convP 2 \eta_2.
\end{align}
We finish the proof by applying Slutsky's theorem.

\subsection{Proof of Corollary~\ref{Corollary: uniform null} and \ref{Corollary: Piecewise uniform alternative}}
These results are direct applications of Theorem~\ref{Thm: Poisson approximation}; hence we omit the proof.

\subsection{Variance of $U_A$} \label{Sec: Variance of U_A}

In the next lemma, we calculate the closed-form of the variance of $U_A$.

\begin{lemma}[Variance of $U_A$] \label{Lemma: Variance}
	Let $A$ be a symmetric positive definite matrix and $\Sigma$ be $\emph{\text{diag}}(\pi) - \pi \pi^{\top}$. Then
	\begin{align}
	\text{\emph{Var}}\left( U_A \right) = \binom{n}{2}^{-1} \Big\{ \emph{tr} \{ (A \Sigma)^2 \}+ 2(n-1) (\pi-\pi_0)^\top A \Sigma A(\pi - \pi_0) \Big\}. \label{Eq: Variance of U_A}
	\end{align}
	Therefore, 
	\begin{align*}
	& \text{\emph{Var}}\left( U_{\pi_0} \right) =\binom{n}{2}^{-1}\Big\{ \emph{tr}\{ (D_{\pi_0} \Sigma)^2 \}+ 2(n-1) (\pi-\pi_0)^\top D_{\pi_0} \Sigma D_{\pi_0} (\pi - \pi_0) \Big\} \quad \text{and} \\[.5em]
	& \text{\emph{Var}}\left( U_I \right) = \binom{n}{2}^{-1} \Big\{ \emph{tr}\left(\Sigma^2 \right)+ 2(n-1) (\pi-\pi_0)^\top \Sigma (\pi - \pi_0) \Big\},
	\end{align*}
	where $D_{\pi_0}$ is defined in (\ref{Eq: chi-square kernel}).
	\begin{proof}
		
		Without loss of generality, we assume $A=I$. Otherwise, define $\bX_1^\ast = A^{1/2} \bX_1$, $\bX_2^\ast = A^{1/2} \bX_2$ and $\pi_0^\ast = A^{1/2} \pi_0$ so that $h_A(\bX_1,\bX_2)  = \left( \bX_1^\ast - \pi_0^\ast  \right) ^\top \left(\bX_2^\ast - \pi_0^\ast \right)$, and proceed similarly. Note that the variance of double summations becomes
		\begin{equation}
		\begin{aligned}
		& \text{{Var}}\left( \sum_{1 \leq i<j \leq n} h_I(\bX_i,\bX_j) \right) \\[.5em] \label{Eq: U-Variance}
		& = \binom{n}{2} \text{Var} \left( h_I(\bX_1, \bX_2) \right) +2(n-2) \text{Cov}\left( h_I(\bX_1, \bX_2), h_I(\bX_1, \bX_3) \right).
		\end{aligned}
		\end{equation}
		We treat the variance and the covariance separately. First, calculate the variance of the kernel:		
		\begin{align*}
		\text{Var} \left( h_I(\bX_1, \bX_2) \right) &  = \mE \left[ h_I^2(\bX_1, \bX_2) \right] - \big\{  \mE \left[ h_I(\bX_1, \bX_2) \right] \big\}^2 \\[.5em]
		&= \mE \left[ h_I^2(\bX_1,\bX_2) \right] - ||\pi - \pi_0||^2_2.
		\end{align*}
		The expected value can be decomposed into:
		\begin{align*}
		\mE \left[ h_I^2(\bX_1,\bX_2) \right]  &= \mE \left[ \big\{ \left(\bX_1 -\pi_0 \right)^\top \left(\bX_2 - \pi_0 \right) \big\}^2 \right] \\[.5em]
		& = \mE \left[ \big\{ \left(\bX_1 - \pi + \pi - \pi_0 \right)^\top \left(\bX_2 - \pi + \pi - \pi_0 \right) \big\}^2 \right] \\[.5em]
		& = \mE \left[ \big\{ (\bX_1-\pi)^\top (\bX_2 - \pi) \big\}^2\right] +  2\mE \left[ \big\{ (\bX_1-\pi)^\top (\pi - \pi_0) \big\}^2 \right] + ||\pi - \pi_0||^2_2.
		\end{align*}
		The first term can be simplified as
		\begin{align*}
		\mE \left[ \big\{ (\bX_1-\pi)^\top (\bX_2 - \pi) \big\}^2\right] & = \mE \left[ \text{tr} \big\{ (\bX_1-\pi)^\top (\bX_2 - \pi) (\bX_2-\pi)^\top (\bX_1 - \pi) \big\} \right] \\[.5em]
		& = \mE \left[ \text{tr} \big\{  (\bX_1 - \pi) (\bX_1-\pi)^\top (\bX_2 - \pi) (\bX_2-\pi)^\top\big\} \right]  \\[.5em]
		& = \text{tr} \Big\{  \mE \left[  (\bX_1 - \pi) (\bX_1-\pi)^\top (\bX_2 - \pi) (\bX_2-\pi)^\top \right] \Big\} \\[.5em]
		& = \text{tr} \Big\{  \mE \left[  (\bX_1 - \pi) (\bX_1-\pi)^\top \right] \mE \left[ (\bX_2 - \pi) (\bX_2-\pi)^\top \right] \Big\}  \\[.5em]
		& = \text{tr}\left( \Sigma^2 \right).
		\end{align*}
		On the other hand, the second term becomes
		\begin{align*}
		\mE \left[ \big\{ (\bX_1-\pi)^\top (\pi - \pi_0) \big\}^2 \right] &= \mE \left[ (\pi - \pi_0) ^\top  (\bX_1-\pi) (\bX_1-\pi)^\top (\pi - \pi_0) \right] \\[.5em]
		& = (\pi - \pi_0) ^\top  \Sigma (\pi - \pi_0). 
		\end{align*}
		Hence, the variance of the kernel can be calculated by
		\begin{align} \label{Eq: Variance Simple Form}
		\text{Var} \left( h_I(\bX_1, \bX_2) \right) = \text{tr}\left( \Sigma^2 \right) + 2 (\pi-\pi_0)^\top \Sigma (\pi-\pi_0).
		\end{align}
		
		\vskip .5em
		
		Next, turn our attention to the covariance:
		\begin{align} \nonumber
		\text{Cov}\left( h_I(\bX_1,\bX_2), h_I(\bX_1,\bX_3) \right) & = \mE \left[ h_I(\bX_1,\bX_2)h_I(\bX_1,\bX_3)  \right] - \mE \left[ h_I(\bX_1,\bX_2)\right] \mE \left[ h_I(\bX_1,\bX_3)\right] \\[.5em]  
		& = \mE \left[ h_I(\bX_1,\bX_2)h_I(\bX_1,\bX_3)  \right] - ||\pi - \pi_0||_2^4. \label{Eq: Cov1}
		\end{align}
		Note that 
		\begin{align} \nonumber
		\mE \left[ h_I(\bX_1,\bX_2)h_I(\bX_1,\bX_3)  \right]  & =  \mE \left[  (\bX_1-\pi_0)^\top (\bX_2 - \pi_0) (\bX_1-\pi_0)^\top (\bX_3 - \pi_0) \right] \\[.5em]
		& =  \mE \left[  (A_1 + B)^\top (A_2 + B) (A_1 + B )^\top (A_3 + B) \right], \label{Eq: covaraince}
		\end{align}
		where $A_i = \bX_i - \pi$ for $i=1,2,3$, and $B = \pi - \pi_0$. In fact, (\ref{Eq: covaraince}) is equivalent to
		\begin{align} \label{Eq: Cov2}
		\mE\left[ A_1^\top B A_1^\top B \right] + B^\top B B^\top B  = (\pi - \pi_0)^\top \Sigma (\pi - \pi_0) + || \pi - \pi_0 ||^4_2,
		\end{align}
		due to 
		\begin{align*}
		\mE\left[ A_1^\top A_2 A_3^\top B \right] & = \mE \left[ \left(\bX_1 - \pi \right)^\top \left(\bX_2 - \pi \right)  \left(\bX_3 - \pi \right)^\top (\pi - \pi_0)  \right] \\[.5em] 
		& =  \mE \left[ \left(\bX_1 - \pi \right)^\top \Big\{ \mE \big[ \left(\bX_2 - \pi \right)  \left(\bX_3 - \pi \right)^\top \big| \bX_1 \big] \Big\}  (\pi - \pi_0)  \right]  \\[.5em]
		& =  \mE \left[ \left(\bX_1 - \pi \right)^\top \Big\{ \mE \left[ \left(\bX_2 - \pi \right) \big| \bX_1 \right]  \mE \big[  \left(\bX_3 - \pi \right)^\top \big| \bX_1 \big] \Big\}  (\pi - \pi_0)  \right]  \\
		& = 0.
		\end{align*}
		In the same way, we can see the other terms become zero. Then, we get a simple form of the covariance from (\ref{Eq: Cov1}) and (\ref{Eq: Cov2}):
		\begin{align}
		\text{Cov}\left( h_I(\bX_1,\bX_2), h_I(\bX_1,\bX_3) \right) = (\pi - \pi_0)^\top \Sigma (\pi - \pi_0). \label{Eq: Covariance Simiple Form}
		\end{align}
		We finish the proof by multiplying $\binom{n}{2}^{-2}$ to (\ref{Eq: U-Variance}) together with (\ref{Eq: Variance Simple Form}) and (\ref{Eq: Covariance Simiple Form}).
	\end{proof}
\end{lemma}

\vskip .8em

\subsection{Proof of Theorem~\ref{Thm: Global Gaussian}}
The proof is based on Corollary 3.1 of \cite{hall1980martingale}. Under the null, $U_A$ is a degenerate centered $U$-statistic, which satisfies $\mE\left[ h_A(\bX_1, \bX_2) \right] = 0$ and $\mE\left[ h_A(\bX_1, \bX_2) \big| \bX_2 \right]=0$. We follow the similar proof steps in Theorem 1 of \cite{hall1984central}, but we adapt the argument to obtain the convergence result for the uniform null in Corollary \ref{Corollary: Uniform Null Global Gaussian}.

\vskip .5em

First, we define the filtration  $\mathcal{F}_k = \sigma\left(\bX_1, \ldots, \bX_k \right)$, and let 
\begin{align*}
Y_j = \sum_{i=1}^{j-1} h_A(\bX_i, \bX_j)  \quad \text{and} \quad  \mS_k =  \sum_{j=2}^k Y_j,
\end{align*}
for $2 \leq k \leq n$. It is easy to check that $\{ \left(\mS_k, \mathcal{F}_k \right) \}$ is a square integrable martingale sequence with zero mean as 
\begin{align*}
\mE \left[ \mS_j \right] = 0 \quad \text{and} \quad \mE\left[ \mS_{i} \big| \mathcal{F}_j \right]  = \mS_j + \sum_{k=j+1}^i \mE \left[  Y_k \big| \mathcal{F}_j   \right] = \mS_j
\end{align*}
for any $i \geq j$. Denote the variance of $\sum_{i<j} h_A(\bX_i, \bX_j)$ by $s_n^2 = {\binom{n}{2}} \text{tr} \{ (A \Sigma)^2 \}$. Then, according to Corollary 3.1 of \cite{hall1980martingale}, it is enough to show that the following two conditions are satisfied under the given assumptions:
\begin{enumerate}\setlength{\itemindent}{.2in}
	\item [\textbf{(C.1)}] \ $\displaystyle s_n^{-2} \sum_{i=2}^n \mE \left[ Y_i^2 I \left( |Y_i | > \epsilon s_n \right) \right] ~\rightarrow~ 0$.
	\item [\textbf{(C.2)}] \ $\displaystyle s_n^{-2} \sum_{i=2}^n \mE \left[ Y_i^2 \big| \mathcal{F}_{i-1}\right] \convP 1$.
\end{enumerate}
Let us first verify the first condition (C.1). Since $|Y_i| > \epsilon s_n$ implies
\begin{align*}
Y_i^2 = \frac{|Y_i|^{2 + \delta}}{|Y_i|^\delta} \leq \frac{|Y_i|^{2+\delta} }{(\epsilon s_n)^\delta},
\end{align*}
for any $\epsilon, \delta >0$, we have
\begin{align*}
s_n^{-2} \sum_{i=2}^n \mE \left[ Y_i^2 I \left( |Y_i | > \epsilon s_n \right) \right]   \leq  \epsilon^{-\delta} s_n^{-2 - \delta} \sum_{i=2}^n {\mE \big[ |Y_i|^{2 + \delta}  \big] }.
\end{align*}
By choosing $\delta = 2$, we will show that $s_n^{-4} \sum_{i=2}^n {\mE \big[ Y_i^{4}  \big] } \rightarrow 0$ to verify (C.1). From the fact that, for any distinct $(i_1, i_2, i_3, i_4)$ or for any combination $(i_1,i_2,i_3,i_4)$ where only one of them is different,
\begin{align*}
\mE \left[ h_A(\bX, \bX_{i_1}) h_A(\bX, \bX_{i_2})h_A(\bX, \bX_{i_3})h_A(\bX, \bX_{i_4})  \right] = 0,
\end{align*}
we can see that
\begin{align*}
\mE \left[ Y_i^4 \right]  & = 	\sum_{i_1, i_2, i_3, i_4 = 1}^{i-1} \mE  \left[ h_A(\bX_i, \bX_{i_1}) h_A(\bX_i, \bX_{i_2}) h_A(\bX_i, \bX_{i_3}) h_A(\bX_i, \bX_{i_4}) \right] \\[.5em]
& = (i-1) \mE\left[  \big\{ h_A(\bX_1, \bX_2) \big\}^4 \right] + 3 (i-1)(i-2) \mE \left[ \big\{ h_A(\bX_1, \bX_2) \big\}^2 \big\{ h_A(\bX_1, \bX_3) \big\}^2 \right].
\end{align*}
Hence, we have
\begin{align*}
\sum_{i=2}^n {\mE \big[ Y_i^{4}  \big] }  & =  \frac{n(n-1)}{2} \mE\left[  \big\{ h_A(\bX_1, \bX_2) \big\}^4 \right] +  n(n-1)(n-2)\mE \left[ \big\{ h_A(\bX_1, \bX_2) \big\}^2 \big\{ h_A(\bX_1, \bX_3) \big\}^2 \right].
\end{align*}
From the second assumption in (\ref{Eq: Gaussian Martingale Assumptions}), it is easy to see $s_n^{-4} \sum_{i=2}^n {\mE \big[ Y_i^{4}  \big] } \rightarrow 0$, which verifies (C.1). 

\vskip .5em

Now, we prove that (C.2) holds under the given conditions, that is to show
\begin{align*}
\frac{2}{n(n-1)\text{tr} \{ (A \Sigma)^2 \}} \sum_{i=2}^n \mE \left[ Y_i^2  \big| \mathcal{F}_{i-1}\right] \convP 1.
\end{align*}
First, we can see from $\mE \left[h_A(\bX_1, \bX_2) h_A(\bX_1, \bX_3) \right] =0$ and $\mE \left[ h_A^2(\bX_1,\bX_2) \right] =\text{tr} \{ (A \Sigma)^2\}$, that
\begin{align*}
\sum_{i=2}^n \mE \left[ Y_i^2  \right] & = \sum_{i=2}^n \sum_{j_1, j_2 =1}^{i-1} \mE \left[ h_A(\bX_i, \bX_{j_1}) h_A(\bX_i, \bX_{j_2}) \right] \\[.5em]
& = \sum_{i=2}^n (i-1) \mE \left[ h_A(\bX_i, \bX_{1}) \right] = \frac{n(n-1)}{2} \text{tr} \{ (A \Sigma)^2 \}.
\end{align*}
Therefore, it is sufficient to prove 
\begin{align*}
\frac{4}{n^2(n-1)^2 \text{tr} \{ (A \Sigma)^2 \}^2} \sum_{i_1, i_2 =2}^n \text{Cov} \left( \mE\left[Y_{i_1}^2 \big| \mathcal{F}_{{i_1}-1} \right], \mE\left[Y_{i_2}^2 \big| \mathcal{F}_{i_2-1}\right] \right) ~\rightarrow~ 0.
\end{align*}
Let us define 
\begin{align} \nonumber
G_A(\bX_i, \bX_j) & = \mE \left[ h_A(\bX_i, \bX_k) h_A(\bX_j, \bX_k) \big| \sigma(\bX_i, \bX_j) \right] \\[.5em]
& = (\bX_i - \pi_0)^\top A \Sigma_0 A (\bX_j -\pi_0), \label{Eq: simple GA}
\end{align}
so that
\begin{align*}
\mE \left[ Y_i^2 | \mathcal{F}_{i-1} \right] & = \sum_{j_1, j_2=1}^{i-1}  \mE \left[ h_A(\bX_i, \bX_{j_1}) h_A(\bX_i, \bX_{j_2}) \big| \mathcal{F}_{i-1} \right] = \sum_{j_1, j_2=1}^{i-1}  G_A(\bX_{j_1}, \bX_{j_2}).
\end{align*}
Then the covariance becomes
\begin{align*}
\text{Cov} \left( \mE\left[Y_{i_1}^2 \big| \mathcal{F}_{{i_1}-1} \right], \mE\left[Y_{i_2}^2 \big| \mathcal{F}_{i_2-1}\right] \right) = \sum_{j_1, j_2=1}^{i_1-1} \sum_{j_1^\prime, j_2^\prime=1}^{i_2-1} \text{Cov} \left(G_A(\bX_{j_1}, \bX_{j_2}), G_A(\bX_{j^\prime_1}, \bX_{j^\prime_2}) \right).  
\end{align*}
Note that for $j_1 \leq j_2$ and $j_1^\prime \leq j_2^\prime$, 
\begin{align*}
\text{Cov} \left(G_A(\bX_{j_1}, \bX_{j_2}), G_A(\bX_{j^\prime_1}, \bX_{j^\prime_2}) \right) =
\begin{cases}
~ \text{Var}\left( G_A(\bX_1, \bX_1) \right) & \text{if}  \ j_1 = j_2 = j_1^\prime = j_2^\prime \\
~\mE \left[ G_A(\bX_1, \bX_2)^2  \right] 	 & \text{if} \ j_1 = j_1^\prime \neq j_2 = j_2^\prime \\
~0 	 & \text{otherwise}
\end{cases}.
\end{align*}
Hence, if $i_1 \geq i_2$, 
\begin{align*}
& \text{Cov} \left( \mE\left[Y_{i_1}^2 \big| \mathcal{F}_{{i_1}-1} \right], \mE\left[Y_{i_2}^2 \big| \mathcal{F}_{i_2-1}\right] \right) =  (i_2-1) \text{Var}\left( G_A(\bX_1, \bX_1) \right) + 2(i_2-1)(i_2-2) \mE \left[ G_A(\bX_1, \bX_2)^2  \right] 
\end{align*}
and the sum of the covariance becomes
\begin{align*}
\sum_{i_1, i_2 =2}^n \text{Cov} \left( \mE\left[Y_{i_1}^2 \big| \mathcal{F}_{{i_1}-1} \right], \mE\left[Y_{i_2}^2 \big| \mathcal{F}_{i_2-1}\right] \right) \leq C_1 \big\{ n^3 \text{Var}\left( G_A(\bX_1, \bX_1) \right) + n^4 \mE \left[ G_A(\bX_1, \bX_2)^2  \right]  \big\},
\end{align*}
where $C_1$ is a constant independent on $n$. Using (\ref{Eq: simple GA}), we have 
\begin{align*}
& \mE \left[ G_A(\bX_1, \bX_2)^2  \right]  = \mE \left[ \big\{  (\bX_1- \pi_0)^\top A \Sigma_0 A (\bX_2 -\pi_0) \big\}^2  \right] = \text{tr} \big\{ ( A\Sigma_0)^4 \big\}, \\[.5em]
& \text{Var}\left( G_A(\bX_1, \bX_1) \right) = \mE \left[ \big\{ h_A(\bX_1, \bX_2) \big\}^2 \big\{ h_A(\bX_1, \bX_3) \big\}^2 \right] - \Big\{\mE \left[   \big\{ h_A(\bX_1, \bX_2) \big\}^2  \right] \Big\}^2 \\[.5em]
&  ~~~~~~~~~~~~~~~~~~~~~~~\leq \mE \left[ \big\{ h_A(\bX_1, \bX_2) \big\}^2 \big\{ h_A(\bX_1, \bX_3) \big\}^2 \right].
\end{align*}
Now, under the given conditions, we bound
\begin{align*}
& \frac{4}{n^2(n-1)^2 \text{tr} \{ (A \Sigma)^2 \}^2} \sum_{i_1, i_2 =2}^n \text{Cov} \left( \mE\left[Y_{i_1}^2 \big| \mathcal{F}_{{i_1}-1} \right], \mE\left[Y_{i_2}^2 \big| \mathcal{F}_{i_2-1}\right] \right)  \\[.5em]
\leq ~ &  C_2 \left( \frac{\text{tr} \{(A\Sigma)^4 \} + n^{-1}\mE \big[ \big\{  h_A(\bX_1, \bX_2)  \big\}^2 \big\{ h_A(\bX_1, \bX_3)  \big\}^2  \big] }{ \text{tr} \{ (A \Sigma)^2 \}^2} \right) ~\rightarrow~ 0
\end{align*}
where $C_2$ is a constant independent on $n$. This completes the proof.

\subsection{Proof of Corollary~\ref{Corollary: Uniform Null Global Gaussian}}
Note that the variance of $U_I$ of the uniform null distribution is
\begin{align*}
{\binom{n}{2}}^{-1} \text{tr} (\Sigma^2) ={\binom{n}{2}}^{-1} \frac{1}{d} \left(1 - \frac{1}{d}\right).
\end{align*}
Therefore, it is enough to show that if $n/\sqrt{d} \rightarrow \infty$, then the conditions of Theorem \ref{Thm: Global Gaussian} are satisfied. To check the first condition, we calculate tr$(\Sigma^4)$ and tr$(\Sigma^2)$ as
\begin{align*}
\text{tr} (\Sigma^4) = \frac{1}{d} \left(1 - \frac{1}{d}\right) \bigg\{ \frac{1}{d^2} \left(1 -\frac{1}{d} \right) + \frac{1}{d^3} \bigg\} \quad \text{and} \quad \text{tr} (\Sigma^2) = \frac{1}{d}\left(1 - \frac{1}{d} \right),
\end{align*}
so that 
\begin{align*}
\frac{\text{tr} (\Sigma^4) }{  \{\text{tr} (\Sigma^2) \}^2} = \frac{1}{d} + \frac{1}{d(d-1)} \rightarrow 0 \quad \text{as } d\rightarrow \infty.
\end{align*}
Next, we verify the second condition when $n / \sqrt{d} \rightarrow \infty$:
\begin{align*}
\frac{ \mE \big[ \big\{ h_I(\bX_1, \bX_2) \big\}^4 \big] + n\mE \big[ \big\{  h_I(\bX_1, \bX_2)  \big\}^2 \big\{ h_I(\bX_1, \bX_3)  \big\}^2  \big]}{ n^2 \{\text{tr} (\Sigma^2) \}^2 } \rightarrow 0.
\end{align*}
For the first part, we have
\begin{align*}
\mE \big[ \big\{ h_I(\bX_1, \bX_2) \big\}^4 \big] & = \mE \big[ \big\{ (\bX_1 - \pi_0)^\top (\bX_2 - \pi_0) \big\}^4 \big] \\[.5em] 
& = \frac{1}{d} \left( 1 - \frac{1}{d} \right)  \Bigg\{ \frac{1}{d^3} + \left(1- \frac{1}{d} \right)^3 \Bigg\}.
\end{align*}
Therefore, 
\begin{align*}
\frac{ \mE \big[ \big\{ h_I(\bX_1, \bX_2) \big\}^4 \big]}{n^2 \{\text{tr} (\Sigma^2) \}^2 } = \frac{\frac{1}{d^3}}{n^2 \frac{1}{d}\left(1 - \frac{1}{d}\right)} + \frac{\left( 1 - \frac{1}{d} \right)^3}{n^2 \frac{1}{d}\left(1 - \frac{1}{d}\right)}  \leq \frac{1}{n^2 d(d-1)} + \frac{d}{n^2} ~\rightarrow~ 0.
\end{align*}
For the second part, 
\begin{align*}
\mE \big[ \big\{  h_I(\bX_1, \bX_2)  \big\}^2 \big\{ h_I(\bX_1, \bX_3)  \big\}^2  \big]  &  = \mE \left[ \big\{ (\bX_1 - \pi_0)^\top (\bX_2 - \pi_0) \big\}^2 \big\{ (\bX_1 - \pi_0)^\top (\bX_3 - \pi_0) \big\}^2 \right] \\[.5em]
&  = \frac{1}{d^2} \left( 1 - \frac{1}{d} \right)^2.
\end{align*}
This gives the second condition:
\begin{align*}
\frac{\mE \big[ \big\{  h_I(\bX_1, \bX_2)  \big\}^2 \big\{ h_I(\bX_1, \bX_3)  \big\}^2  \big]}{ n \{\text{tr} (\Sigma^2) \}^2 } = \frac{1}{n} ~\rightarrow~ 0.
\end{align*}
Hence, the proof is complete.

\subsection{Proof of Theorem~\ref{Thm: Gaussian Limit under alternative}}

Note that the explicit formula for $\text{Var}(U_A)$ is established in Lemma~\ref{Lemma: Variance}. Recall the decomposition $U_A = U_{\text{quad}} + U_{\text{linear}}$ given in the main text. Then under $(S.1)$, we have 
\begin{align*}
\frac{U_A - || A^{1/2}(\pi - \pi_0) ||_2^2  }{\sqrt{\text{Var}(U_A) }} = \frac{ U_{\text{linear}} - || A^{1/2}(\pi - \pi_0) ||_2^2  }{ \sqrt{\text{Var} (U_{\text{linear}})}} + o_P(1),
\end{align*}
and the asymptotic normality follows by the usual central limit theorem. On the other hand, under (S.2), we have
\begin{align*}
\frac{U_A - || A^{1/2}(\pi - \pi_0) ||_2^2  }{\sqrt{\text{Var}(U_A) }} = \frac{ U_{\text{quad}} }{ \sqrt{\text{Var} (U_{\text{quad}})}} + o_P(1).
\end{align*}
Then we follow the similar steps in the proof of Theorem \ref{Thm: Global Gaussian} to get the normality of $U_{\text{quad}}$. Hence the proof is complete.

%\subsection{Proof of Theorem~\ref{Thm: Global Minimax Optimality based on a mixture distribution}}
%This result is a direct consequence of Theorem~\ref{Thm: Global minimax optimality of $U_{w}$}; hence omitted.

\subsection{Proof of Theorem~\ref{Thm: Global minimax optimality of $U_{w}$}}
We proceed along the lines of the proof of Theorem 2 in \cite{balakrishnan2017hypothesis}. Note that the expectation and variance of $U_{w}$ (Lemma~\ref{Lemma: Variance}) are given by
\begin{align*}
& \mE\left[U_{w} \right] = || A_{w}^{1/2} (\pi - \pi_0) ||_2^2 \\[.5em]
& \text{Var}\left[ U_{w} \right] = {\binom{n}{2}}^{-1} \Big\{  \text{tr} \{ (A_{w} \Sigma)^2 \} + 2(n-1) (\pi-\pi_0)^\top A_{w} \Sigma A_{w}(\pi - \pi_0)  \Big\}.
\end{align*}
Let $\mE_0[\cdot]$, $\mE_1[\cdot]$ be the expected value under the null and the alternative, respectively, and similarly denote Var$_0[\cdot]$, Var$_1[\cdot]$. By Chebyshev's inequality, under the null, we can see
\begin{align*}
\mP_{H_0} \left( U_{w}  \geq t_\alpha \right) \leq \frac{\text{Var}_0\left[U_{w} \right]}{t_\alpha^2} = \alpha, 
\end{align*}
and $t_\alpha = \sqrt{\alpha^{-1} \text{Var}_0\left[U_{w} \right]}$. This shows that $\phi(U_{w})$ has size at most $\alpha$. 

\vskip .5em

For the type II error bound, assume the following two conditions are true:
\begin{align*}
(i) \ \ t_\alpha  \leq \frac{\mE_1 \left[ U_{w} \right] }{2} \quad \text{and} \quad (ii) \ \  \sqrt{\frac{\text{Var}_1 \left[ U_{w}\right]}{\zeta}} \leq \frac{\mE_1 \left[ U_{w}\right]}{2}.
\end{align*} 
Then, we can observe that 
\begin{align*}
\mP_{H_1} \left(\phi(U_{w}) = 0\right) & ~ = ~ \mP_{H_1} \left(U_{w} < t_\alpha \right) \\[.5em] 
& ~ \leq ~  \mP_{H_1} \left(U_{w} < \frac{\mE_1 \left[ U_{w} \right] }{2} \right) & \quad \text{by (\emph{i})} \\[.5em]
& ~ = ~ \mP_{H_1} \left(U_{w} < \mE_1 \left[ U_{w} \right]  - \frac{\mE_1 \left[ U_{w} \right] }{2} \right)  \\[.5em]
& ~ \leq ~ \mP_{H_1} \left(U_{w} < \mE_1 \left[ U_{w} \right]  - \sqrt{\frac{\text{Var}_1 \left[ U_{w}\right]}{\zeta}}\right) & \quad \text{by (\emph{ii})} \\[.5em]
& ~ \leq ~ \zeta,
\end{align*}
where the last inequality follows by Chebyshev's inequality. Therefore, the proof can be done by showing $(i)$ and $(ii)$. 

\vskip .5em

We begin with proving the first part $(i)$. After some calculations, we can see
\begin{align} \label{Eq: Null Variance Simple Form}
& \text{tr} \{ (A_w \Sigma)^2 \} ~=~ \sum_{j=1}^d \frac{\pi_j^2}{w_j^2} -2 \sum_{j=1}^d \frac{\pi_j^3}{w_j^2} + \left(\sum_{j=1}^d \frac{\pi_j^2}{w_j}\right)^2.
\end{align}
Therefore, under the null, the variance of $U_w$ can be expanded to
\begin{align*}
\text{Var}_0\left[ U_{w} \right] = {\binom{n}{2}}^{-1}\Bigg\{ \sum_{j=1}^d \frac{\pi_{0,j}^2}{w_j^2} -2 \sum_{j=1}^d \frac{\pi_{0,j}^3}{w_j^2} + \left(\sum_{j=1}^d \frac{\pi_{0,j}^2}{w_j}\right)^2    \Bigg\}.
\end{align*}
By Cauchy-Schwarz inequality, note that
\begin{align} \label{Eq: Upper bound of trace term}
\left(\sum_{j=1}^d \frac{\pi_j^{2}}{w_j} \right)^2  =~ \left(\sum_{j=1}^d \frac{\pi_j^{3/2}}{w_j}  \pi_j^{1/2}\right)^2  ~ \leq ~ \sum_{j=1}^d \frac{\pi_j^3}{w_j^2} \sum_{j=1}^d \pi_j ~=~ \sum_{j=1}^d \frac{\pi_j^3}{w_j^2},
\end{align}
which implies 
\begin{align} \label{Eq: Upper bound of the variance}
\text{Var}_0\left[ U_{w} \right] ~ \leq ~ {\binom{n}{2}}^{-1} \sum_{j=1}^d \frac{\pi_{0,j}^2}{w_j^2}.
\end{align}
Using (\ref{Eq: Upper bound of the variance}), the critical value is upper bounded by
\begin{align*}
t_\alpha  =  \sqrt{\alpha^{-1} \text{Var}_0 \left[ U_w \right]}  ~ \leq  \frac{2}{n}\sqrt{\frac{1}{\alpha}\sum_{j=1}^d \frac{\pi_{0,j}^2}{w_j^2}}
\end{align*}
Note that, from the comparable condition, there exist $C_1, C_2 >0$ and $\gamma \in (0,1)$ such that 
\begin{align} \label{Eq: Comparable Condition}
C_1 \{  \gamma \pi_{0,i} + (1-\gamma) 1/d \} \leq w_i \leq  C_2 \{ \gamma \pi_{0,i} + (1-\gamma) 1/d\} \quad \text{for all $i = 1,\ldots,d$.}
\end{align}
Consequently, the critical value is further upper bounded by 
\begin{align} \nonumber
t_\alpha & \leq ~ \frac{2}{n} \sqrt{\frac{1}{\alpha}\sum_{j=1}^d \frac{\pi_{0,j}^2}{w_j^2}} ~ \leq ~  \frac{2}{n} \sqrt{\frac{1}{\alpha}\sum_{j=1}^d \left( \frac{\pi_{0,j}}{C_1 \{  \gamma \pi_{0,i} + (1-\gamma) 1/d \}} \right)^2} \\[.5em] \label{Eq: upper bound of talp}
& \leq ~ \frac{2}{C_1 \gamma n} \sqrt{ \frac{d}{\alpha}},
\end{align}
where the last inequality is due to $1/ \{  C_1 \gamma \pi_{0,i} + C_1 (1-\gamma) 1/d \} \leq 1 / \{ C_1 \gamma \pi_{0,i} \}$. On the other hand, Cauchy-Schwarz inequality together with the comparable condition presents
\begin{align} \label{Eq: lower bound of expected value}
\mE_1 \left[ U_w \right] =  \sum_{i=1}^d \frac{(\pi_{i} - \pi_{0,i})^2}{w_i} ~ \geq ~ \frac{|| \pi - \pi_0 ||_1^2}{\sum_{i=1}^d w_i} ~ \geq ~ \frac{\epsilon_n^2}{C_2}.
\end{align}
Therefore, the first condition $(i)$ is satisfied if
\begin{align*}
\epsilon_n^2~ \geq ~ \frac{4 C_2}{C_1 \gamma n} \sqrt{ \frac{d}{\alpha}}.
\end{align*}
This is the case from the assumption in (\ref{Eq: Condition for critical radius}).

\vskip .8em

Next, we prove the condition $(ii)$. First, observe that
\begin{align*}
\text{Var}_1 \left[ U_{w} \right] = \frac{1}{\binom{n}{2}} \Big\{ \text{tr} \{ (A_{w} \Sigma)^2 \}+ 2(n-1) (\pi-\pi_0)^\top A_{w} \Sigma A_{w} (\pi - \pi_0) \Big\}.
\end{align*}
By using the result in (\ref{Eq: Null Variance Simple Form}) and (\ref{Eq: Upper bound of trace term}), the first trace term is bounded by
\begin{align*}
\frac{1}{\binom{n}{2}} \text{tr} \{ (A_{w} \Sigma)^2 \} ~ \leq ~ \frac{4}{n^2}\sum_{j=1}^{d} \frac{\pi_j^2}{w_j^2}, 
\end{align*}
for $n \geq 2$. On the other hand, the second term is bounded by
\begin{align*}
\frac{4}{n} (\pi-\pi_0)^\top A_{w} \Sigma A_{w} (\pi - \pi_0) & ~ = ~  \frac{4}{n} (\pi-\pi_0)^\top A_{w} \left(\text{diag}\{\pi\} - \pi \pi^\top \right) A_{w} (\pi - \pi_0)  \\[.5em]
& ~ \leq ~ \frac{4}{n} (\pi-\pi_0)^\top A_{w} \text{diag}\{\pi\} A_{w} (\pi - \pi_0) = \frac{4}{n} \sum_{j=1}^d \frac{\Delta_j^2 \pi_j }{w_j^2},
\end{align*}
where $\Delta_i = \pi_{0,i} - \pi_i$. Therefore, we have
\begin{align*}
\text{Var}_1 \left[ U_{w} \right] & ~\leq ~ \frac{4}{n^2}\sum_{j=1}^{d} \frac{\pi_j^2}{w_j^2} + \frac{4}{n} \sum_{j=1}^d \frac{\Delta_j^2 \pi_j }{ w_j^2} \\[.5em]
& ~ = ~  \frac{4}{n^2}\sum_{j=1}^{d} \frac{\pi_{0,j}^2 + \Delta_j^2 - 2 \pi_{0,j} \Delta_j  }{w_j^2} + \frac{4}{n} \sum_{j=1}^d \frac{\Delta_j^2 \pi_{0,j} - \Delta_j^3 }{ w_j^2} \\[.5em]
& ~ \leq ~  \underbrace{\frac{8}{n^2}\sum_{j=1}^{d} \frac{\pi_{0,j}^2}{w_j^2}}_{U_1} + \underbrace{\frac{8}{n^2}\sum_{j=1}^{d} \frac{\Delta_j^2}{w_j^2}}_{U_2}  + \underbrace{\frac{8}{n} \sum_{j=1}^d \frac{\Delta_j^2 \pi_{0,j} }{ w_j^2}}_{U_3} + \underbrace{\frac{8}{n} \sum_{j=1}^d \frac{|\Delta_j |^3 }{ w_j^2}}_{U_4}. 
\end{align*}
To finish the proof, we need to verify 
\begin{align*}
\sum_{i=1}^4 \frac{2\sqrt{U_i/\zeta}}{\mE_1 \left[ U_{w}\right]} \leq 1.
\end{align*}
Indeed, this is the case by modifying the result in \cite{balakrishnan2017hypothesis} with a different constant factor. To show the details, using (\ref{Eq: lower bound of expected value}), the first term is upper bounded by
\begin{align*}
\frac{2\sqrt{U_1/\zeta}}{\mE_1 \left[ U_{w} \right]} ~ \leq ~ \frac{4 \sqrt{2} C_2 }{\sqrt{\zeta} n \epsilon_n^2} \sqrt{\sum_{j=1}^{d} \frac{\pi_{0,j}^2}{w_j^2}} ~ \leq ~ \frac{4 \sqrt{2} C_2}{\sqrt{\zeta}C_1 \gamma} \frac{\sqrt{d}}{n \epsilon_n^2} ~ \leq ~ \frac{1}{4}.
\end{align*}
For the second term, note that
\begin{align*}
U_2  & = \frac{8}{n^2} \sum_{j=1}^{d} \frac{\Delta_j^2}{w_j^2}  \leq \frac{8}{n^2} \sum_{j=1}^{d} \frac{\Delta_j^2}{C_1^2 \{ \gamma \pi_{0,j} + (1-\gamma) 1/d \}^2} = \frac{8d^2}{n^2 C_1^2(1-\gamma)^2} \sum_{j=1}^{d} \frac{\Delta_j^2}{\{ d \pi_{0,j} \gamma / (1-\gamma) + 1 \}^2}  \\[.8em]
& \leq \frac{8d^2}{n^2 C_1^2(1-\gamma)^2} \sum_{j=1}^{d} \frac{\Delta_j^2}{\{ d \pi_{0,j} \gamma / (1-\gamma) + 1 \}} = \frac{8d}{n^2 C_1^2 (1- \gamma)} \underbrace{\sum_{j=1}^d \frac{\Delta_j^2}{\gamma \pi_{0,j}  + (1-\gamma)1/d}}_{\overset{let}{=} \rho_n}.
\end{align*}
The expected value is lower bounded in terms of $\rho_n$ by 
\begin{align*}
\mE_1 \left[U_w \right] = \sum_{j=1}^d \frac{\Delta_j^2}{w_j}  \geq \frac{\rho_n}{C_2},
\end{align*}
and similarly to (\ref{Eq: lower bound of expected value}), it is seen that $\rho_n \geq \epsilon_n^2$. 
Using these results, 
\begin{align}
\frac{2\sqrt{U_2/\zeta}}{\mE_1 \left[ U_{w} \right]}  \leq \frac{4 C_2\sqrt{2d}}{C_1 \sqrt{\zeta(1-\gamma)}n\epsilon_n}  \leq \frac{1}{4}.
\end{align}
For the third term, note that 
\begin{align*}
\frac{\pi_{0,j}}{w_j} \leq \frac{\pi_{0,j}}{C_1 \{ \gamma \pi_{0,j} + (1-\gamma) 1/d  \} } \leq \frac{1}{C_1 \gamma}.
\end{align*}
Using this inequality, 
\begin{align*}
U_3  = \frac{8}{n} \sum_{j=1}^d \frac{\Delta_j^2 \pi_{0,j} }{ w_j^2} \leq \frac{8}{C_1 \gamma n} \sum_{j=1}^d \frac{\Delta_j^2  }{ w_j} = \frac{8}{C_1 \gamma n}\mE_1 \left[ U_w \right].
\end{align*}
As a result, 
\begin{align*}
\frac{2\sqrt{U_3/\zeta}}{\mE_1 \left[ U_{w} \right]}  \leq  \frac{2\sqrt{2 C_2}}{\sqrt{C_1 \gamma \zeta n} \epsilon_n} \leq \frac{1}{4}.
\end{align*}
To control the last term, the monotonicity of the $\ell_p$ norms and the comparable condition present
\begin{align*}
U_4 & = \frac{8}{n} \sum_{j=1}^d \frac{|\Delta_j |^3 }{ w_j^2} \leq \frac{8}{n} \left( \sum_{j=1}^d  \frac{ \Delta_j^2 }{ w_j^{4/3}} \right)^{3/2} \leq ~\frac{8d^{1/2}}{nC_1^{2} (1 - \gamma)^{1/2}}   \left( \sum_{j=1}^d \frac{\Delta_j^2}{\gamma \pi_{0,j} + (1 - \gamma) 1/d}  \right)^{3/2}.
\end{align*}
Based on the result, 
\begin{align*}
\frac{2\sqrt{U_4/\zeta}}{\mE_1 \left[ U_{w} \right]} ~ \leq ~ \frac{4\sqrt{2} C_2 d^{1/4}}{C_1 (1 - \gamma)^{1/4} \sqrt{\zeta n} \rho_n^{1/4}}  ~ \leq ~ \frac{4\sqrt{2} C_2 d^{1/4}}{C_1 (1 - \gamma)^{1/4} \sqrt{\zeta n} \epsilon_n^{1/2}}  \leq \frac{1}{4}. 
\end{align*}
This completes the proof.

\end{document}